# On the uniform distribution of rational inputs with respect to condition numbers of Numerical Analysis[1]

D. Castro[2], J.L. Montaña[2], L. M. Pardo[2], J. San Martín[2]

November 19, 2018

## Abstract

We show that rational data of bounded input length are uniformly distributed with respect to condition numbers of numerical analysis. We deal both with condition numbers of Linear Algebra and with condition numbers for systems of multivariate polynomial equations. For instance, we show that for any $w > 1$ and for any $n \times n$ rational matrix $M$ of bit length $O(n^4 \log n) + \log w$, the condition number $k(M)$ satisfies $k(M) \leq w n^{5/2}$ with probability at least $1 - 2w^{-1}$. Similar estimates are shown for the condition number $\mu_{norm}$ of M. Shub and S. Smale when applied to systems of multivariate homogeneous polynomial equations of bounded input length. Finally we apply these techniques to show the probability distribution of the precision (number of bits of the denominator) required to write down approximate zeros of affine systems of multivariate polynomial equations of bounded input length.

---

[1]: Research was partially supported by the Spanish grants : BFM2000-0349, HF1999-055 and "Acción Integrada Hispano–Argentina Universidad de Cantabria–Universidad de Buenos Aires".

[2]: Departamento de Matemáticas, Estadística y Computación.
Facultad de Ciencias, Universidad de Cantabria, E-39071 Santander, Spain

# Contents

# 1 Introduction

In [4] we exhibited examples of systems of multivariate polynomial equations with integer coefficients of small input length such that the required precision to write down an *approximate zero* (in the sense of [33]) is exponential in the number of variables (cf. also Section 6 below). These examples are based on the fact that the condition number $\mu_{norm}$ of [32] and [3] is doubly exponential in the number of variables.

In his seminal paper of 1948 (cf. [36]), A. Turing introduced two condition numbers in Linear Algebra : the condition numbers $M$ and $N$. A. Turing also showed that the $M-$condition number is a sharp worst case estimate of the loss of precision in linear Algebra Computations. The same condition number was also discussed by J. von Neumann and collaborators (cf. [22]) and by J. H. Wilkinson (cf. [39], [38]).

In his paper, A. Turing already exhibited examples of systems of linear equations with rational coefficients of small input length such that the condition number becomes too large.

In all these examples, input length means the number of tape cells (also digits) required to represent the input in a Turing machine.

Nevertheless, there are deep and meaningful studies on the probability distribution of condition numbers of numerical analysis. The reader may follow the estimates for the Linear Algebra case in [34], [7, 8, 9] and the estimates for the Non-linear case in [31], [3] or [32]. Most of these studies show that condition numbers are (almost) polynomial in the dimension of the source space with high probability. Hence, our examples in [4] and Turing's examples seem to be in contradiction with these studies.

This "contradiction" is misleading. The probability distribution obtained in the series of papers [34], [7, 8, 9], [31, 32], [3] is based on a "continuous model for numerical analysis procedures". A "continuous model" assumes that inputs belong to some "continuous source space". Consequently, studies on the probability distribution under a "continuous model" mean "probability distribution of condition numbers of inputs in some continuous source space".

For instance, in the Linear Algebra case, "inputs" are projective points in the "continuous source space" $\mathbb{P}(\mathcal{M}_n(\mathbb{R}))$, where $\mathcal{M}_n(\mathbb{R})$ is the vector space of all $n \times n$ real matrices. Thus, if $k(M)$ denotes the standard condition number of matrices $M \in \mathcal{M}_n(\mathbb{R})$ (see more details below), the "probability distribution" of $k$ means the probability distribution of $k$ when $M$ runs over the real projective space $\mathbb{P}(\mathcal{M}_n(\mathbb{R}))$.

On the other hand, real life computing is "discrete". This essentially means that real life computers take as input matrices $M \in \mathbb{P}(\mathcal{M}_n(K))$ where $K$ is some "discrete field" ($K = \mathbb{Q}$, for instance). Moreover, an essential parameter in real life computing is the input length. The input length of $M \in \mathcal{M}_n(\mathbb{Q})$ is the number of digits required to represent the entries of $M$ in a computer (cf. also Subsection 2.4 below). Hence, real life computing demands the probability distribution of condition numbers of numerical analysis for inputs of bounded length in a "discrete space" as, for instance, the rational projective space $\mathbb{P}(\mathcal{M}_n(\mathbb{Q}))$. These pages are devoted to show this kind of estimates within a discrete context. We deal both with the Linear and the Non–linear cases.

The rest of the Introduction is devoted to state the new outcomes in a succinct form. Subsection 1.1 below shows the probability distribution of condition numbers $k$ and $\mu$ for rational matrices with bounded input length. Subsection 1.2 shows the probability distribution of condition number $\mu_{norm}$ of M. Shub and S. Smale for systems of multivariate homogeneous polynomial equations of bounded input length whose coefficients are Gauss rationals. Finally, in Subsection 1.3 below we show the probability distribution of the precision required to write down affine approximate zeros of systems of polynomial equations of bounded input length.

An advised reader may realize that our statements imply that rational data of bounded input length are uniformly distributed with respect to condition numbers of numerical analysis. In fact, our statements below exhibit concrete estimates for the discrepancy (i.e. for the difference between the "continuous" and the "discrete" probabilities). In other words, we show the uniform distribution of rational data but paying attention to A. Schonhage's Lemma : "Do care about constants!".

## 1.1 The Linear Algebra Case

For every square matrix $M \in \mathcal{M}_n(\mathbb{R})$, let $k(M)$ be the standard condition number of $M$. This condition number is a bound for the loss of precision of most linear algebra procedures. It was introduced at the end of the forties by several authors : A. Turing in [36], J. von Neumann and coauthors (cf. [22]) and J. H. Wilkinson (cf. [39, 38]). It has been extensively studied in [34], [15], [7, 8, 9]. The condition number $k(M)$ is given by the following identity
$$k(M) := \|M^{-1}\| \|M\|,$$
where $\|M\|$ is the norm of $M$ as linear operator, namely
$$\|M\| := sup\{\|Mv\| \ : \ v \in \mathbb{R}^n, \|v\| = 1\}.$$

The condition number $k$ is naturally defined on the projective space $\mathbb{P}(\mathcal{M}_n(\mathbb{R}))$. Thus, we denote by $k$ the mapping
$$k : \mathbb{P}(\mathcal{M}_n(\mathbb{R})) \longrightarrow \mathbb{R}_+.$$

Let $\mathbb{P}(\mathcal{M}_n(\mathbb{Q}))$ be the projective space defined by all $n \times n$ matrices with rational entries.

In these pages, the input length of a projective point $P$ in $\mathbb{P}(\mathcal{M}_n(\mathbb{Q}))$ defined by a matrix $M \in \mathbb{P}(\mathcal{M}_n(\mathbb{Q}))$ is given as the logarithm of the Northcott–Schmidt height of the projective point $P \in \mathbb{P}(\mathcal{M}_n(\mathbb{Q}))$ (cf. Subsection 2.4 below).

The following statement holds :

**Theorem 1 (Main Theorem for Linear Algebra)** *With the same notations and assumptions as above, there is a universal positive constant $c_1 \in \mathbb{R}$, $c_1 \leq 10$ such that the following holds :*
*Let $h, w \in \mathbb{R}$ be positive real numbers. Assume that $w > 1$ and $h$ is big enough to satisfy the following inequality :*

$$h \geq c_1 n^4 \log n + \log w.$$

*Then, for every matrix $M \in \mathbb{P}(\mathcal{M}_n(\mathbb{Q}))$ of bit length at most $h$, the condition number $k(M)$ satisfies $k(M) < wn^{5/2}$ with probability at least*

$$1 - \frac{2}{w}.$$

The same assertion holds for the condition number $\mu$ introduced by S. Smale in [34] (cf. Subsection 2.2 below) and for the loss of precision of the gradient method (cf. [12]).

The following constants are going to be used in what follows. Let $\ell \in \mathbb{N}$ be a positive integer number. By $K_\ell \in \mathbb{R}$ we denote the volume of the $\ell$–dimensional unit ball in $\mathbb{R}^\ell$. Namely,

$$K_\ell := \frac{2\pi^{\frac{\ell}{2}}}{\ell \Gamma\left(\frac{\ell}{2}\right)}.$$

For $\ell = 0$, we write $K_0 := 1$. Finally, let us introduce the constant $\mathfrak{S}^{(m)}$ given by the following identity

$$\mathfrak{S}^{(m)} := \sum_{\ell=0}^{m-1} \binom{m}{\ell} K_\ell. \tag{1}$$

Using E. Artin's estimates for the gamma function (cf. [2]) the following estimates also hold :

$$\mathfrak{S}^{(m)} \leq \frac{1}{\sqrt{\pi}}(1 + \sqrt{2e\pi})^m.$$

and

$$\frac{\mathfrak{S}^{(m)}}{K_m} \leq e^{\frac{2}{3}} m^{\frac{m}{2}} e^{\frac{m}{12}}.$$

Theorem 1 above is a direct consequence of the following more technical (but more complete) statement.

**Theorem 2** *With the same notations as, above let $H, \varepsilon \in \mathbb{R}$ be two positive real numbers. Assume that $H \geq 1$ and $\varepsilon > 0$. For every matrix $M \in \mathbb{P}(\mathcal{M}_n(\mathbb{Q}))$ of Northcott–Schmidt's height at most $H$, the condition number $k(M)$ satisfies $k(M) \geq 1/\varepsilon$ with probability at most*

$$\varepsilon n^{5/2} + \frac{B(\varepsilon, n)}{H - C(n)},$$

*where*

$$B(\varepsilon, n) := \zeta(n^2) \left[ \frac{(T_n + \varepsilon n^{5/2})\mathfrak{S}^{(n^2)}}{K_{n^2}\zeta(n^2 - 1)} + 2\varepsilon n^{5/2} + \frac{2}{K_{n^2}} \right],$$

$$C(n) := \zeta(n^2) \left[ \frac{\mathfrak{S}^{(n^2)}}{K_{n^2}\zeta(n^2 - 1)} + 1 + \frac{1}{K_{n^2}} \right],$$

$\zeta$ *is Riemann's zeta function and*

$$T_n := (2\max\{4, n\})^{(2n^4 + 4n^2)}.$$

## 1.2 The Non–Linear Case

Let $\mathbb{Q}[i]$ be the field of Gauss rationals. For every positive integer number $d \in \mathbb{N}$, let $H_d$ be the complex vector space of all homogeneous polynomials with complex coefficients of degree $d$. Namely,

$$H_d := \{ f \in \mathbb{C}[X_0, \ldots, X_n] : f \text{ homogeneous, } \deg(f) = d \}.$$

For every list of degrees $(d) := (d_1, \ldots, d_n) \in \mathbb{N}^n$, let $\mathcal{H}_{(d)}$ be the complex vector space given by the following identity :

$$\mathcal{H}_{(d)} := H_{d_1} \times \cdots \times H_{d_n}.$$

The space $\mathcal{H}_{(d)}$ is the space of systems of complex homogeneous polynomials $F := (f_1, \ldots, f_n)$ such that $\deg(f_i) = d_i$.
For every list of degrees $(d)$, let $D_{(d)}$ be the maximum $D_{(d)} := \max\{d_1, \ldots, d_n\}$ and let $\mathcal{D}_{(d)}$ be the Bézout number of $(d)$, namely $\mathcal{D}_{(d)} := \prod_{i=1}^n d_i$.
Let $N_d$ be the dimension of $\mathcal{H}_d$. We have

$$N_d := \binom{d+n}{n}.$$

For every system $F \in \mathcal{H}_{(d)}$, let $\mu_{norm}(F)$ be the normalized condition number introduced in [29, 32] (cf. also [3] and Subsection 4.3 below). This is a projective invariant (i.e. $\mu_{norm}(F)$ depends only on the class defined by $F$ in the complex projective space $\mathbb{P}(\mathcal{H}_{(d)})$).

Let $N$ be the (complex) dimension of the projective space $\mathbb{P}(\mathcal{H}_{(d)})$. Observe that the following equality holds :

$$N = \left(\sum_{i=1}^{n} N_{d_i}\right) - 1 = \left(\sum_{i=1}^{n} \binom{d_i + n}{n}\right) - 1.$$

For every list of degrees $(d) := (d_1, \ldots, d_n)$, we denote by $\mathcal{H}_{(d)}(\mathbb{Q}[i])$ the set of systems of homogeneous polynomials $F := (f_1, \ldots, f_n)$ such that $f_i \in \mathbb{Q}[i][X_0, \ldots, X_n]$.

Let $\mathbb{P}(\mathcal{H}_{(d)}(\mathbb{Q}[i]))$ be the projective space defined by the $\mathbb{Q}[i]-$vector space $\mathcal{H}_{(d)}(\mathbb{Q}[i])$. In this case the *(unitarily invariant) bit length of a projective point* $F \in \mathbb{P}(\mathcal{H}_{(d)}(\mathbb{Q}[i]))$ *is defined as the logarithm of the unitarily invariant height of* $F$ (cf. Subsection 4.4 below). Then, the following statement holds :

**Theorem 3 (Main Theorem for the Non–Linear Case)** *With the same notations and assumptions as above, there is a positive universal constant $c_2 \in \mathbb{R}$, $c_2 \leq 20$ such that the following holds :*
*Let $(d)$ be a list of degrees and let $h, w \in \mathbb{R}$ be two positive real numbers. Assume that $w > 1$ and $h$ is big enough to satisfy the following inequality :*

$$h \geq c_2 N^2 \log N + \log w.$$

*Then, for every system of polynomial equations $F \in \mathbb{P}(\mathcal{H}_{(d)}(\mathbb{Q}[i]))$ of (unitarily invariant) bit length at most $h$, the condition number $\mu_{norm}(F)$ satisfies*

$$\mu_{norm}(F) < (w\mathfrak{C}[(d)])^{1/4}, \qquad (2)$$

*with probability at least*

$$1 - \frac{2}{w},$$

*where*

$$\mathfrak{C}[(d)] := n^3(n+1)N(N-1)\prod_{i=1}^{n} d_i.$$

The reader himself may rewrite this statement for the quadratic case (i.e. $(d) = (2, 2, \ldots, 2)$). Then he will realize what we mean.

This statement is an immediate consequence of the following more technical (but more complete) statement. In order to state this new Theorem, let us introduce some more notations :
Let $S_d$ be the set of multi–indices

$$S_d := \{(\mu_0, \ldots, \mu_n) \in \mathbb{N}^{n+1} : |\mu| = \mu_0 + \ldots + \mu_n \leq d\}.$$

We introduce a well–ordering in the set of multi-indices $S_d$ as a bijection

$$\phi : S_d \longrightarrow \{i \in \mathbb{N} : 1 \leq \imath \leq N_d\}.$$

Using this well–ordering on $S_d$, we define the diagonal matrix $\Delta_d$ in the following terms :

$$\Delta_d := \left( \binom{d}{\mu_0 \ldots \mu_n}^{-\frac{1}{2}} \right)_{1 \leq \phi(\mu_0,\ldots,\mu_n) \leq N_d},$$

where

$$\binom{d}{\mu_0 \ldots \mu_n} := \left( \frac{d!}{\mu_0! \cdots \mu_n!} \right).$$

For every list of degrees $(d) := (d_1, \ldots, d_n)$, let $\Delta_{(d)}$ be the diagonal matrix given as the diagonal sum of $\Delta_{d_1}, \ldots, \Delta_{d_n}$. Namely,

$$\Delta_{(d)} := \Delta_{d_1} \oplus \ldots \oplus \Delta_{d_n}. \tag{3}$$

**Theorem 4** *With the same notations as above, let $H, \varepsilon \in \mathbb{R}$ be two positive real numbers. Assume that $H \geq 1$ and $\varepsilon > 0$. Let $(d) := (d_1, \ldots, d_n)$ be a list of degrees.*
*Then, for every system $F \in \mathbb{P}(\mathcal{H}_{(d)}(\mathbb{Q}[i]))$ of unitarily invariant height at most $H$ the condition number $\mu_{norm}(F)$ satisfies*

$$\mu_{norm}(F) \geq 1/\varepsilon$$

*with probability at most*

$$\varepsilon^4 \mathfrak{C}[(d)] + \frac{\mathcal{F}(\varepsilon, N)}{\zeta(2N+2)^{-1} H - \mathcal{G}(N)}$$

*where $\mathfrak{C}[(d)]$ is the constant introduced in Theorem 3 above and*

$$\mathcal{G}(N) := \frac{\mathfrak{S}^{(2N+2)}}{\det(\Delta_{(d)})\zeta(2N+1)K_{2N+2}} + \frac{1}{K_{2N+2}} + 2,$$

$$\mathcal{F}(\varepsilon, N) := \frac{(\overline{T}_N + \varepsilon^4 \mathfrak{C}[(d)])\mathfrak{S}^{(2N+2)}}{\det(\Delta_{(d)})\zeta(2N+1)K_{2N+2}} + \frac{(\varepsilon^4 \mathfrak{C}[(d)] + 1)}{K_{2N+2}} + 4\varepsilon^4 \mathfrak{C}[(d)],$$

*and*

$$\overline{T}_N := \max\{8, 4(D_{(d)} + 1)\}^{4[N^2 + 3(n+2)^2(N+1)]}.$$

### 1.3 An application

We apply the estimates of Theorem 3 above to estimate the probability distribution of the precision required to write down affine approximate zeros in $\mathbb{Q}[i]^n$ of systems of multivariate polynomial equations of bounded input length. The reader may look at Section 6 for detailed definitions. For every $z \in \mathbb{Q}[i]^n$, the precision is defined as the bit length of a denominator of $z$. For every system of polynomial equations $F \in \mathcal{P}_{(d)}(\mathbb{Q}[i])$, the precision $Pr(F)$ of $F$ is defined as the maximum precision required to write down approximate zeros in $\mathbb{Q}[i]^n$ for every actual zero of $F$. Then, we show the following statement.

**Corollary 5** *With the same notations as above, there is a universal constant $c_3 > 0$ ( $c_3 \leq 20$) such that the following holds :*
*Let $(d) := (d_1, \ldots, d_n)$ be a list of degrees and let $h, w \in \mathbb{R}$ be two positive real numbers. Assume that $w > 1$ and $h$ is big enough to satisfy the following inequality :*
$$h \geq c_3 N^2 (\log N + \log D_{(d)}) + \log w,$$

*Then, for every system of polynomial equations $F \in \mathbb{P}(\mathcal{P}_{(d)}(\mathbb{Q}[i]))$ of bit length at most $h$, the required precision $Pr(F)$ to write down an approximate zero of $F$ satisfies*
$$Pr(F) \leq O(n \log_2 D_{(d)} + \log_2 w), \tag{4}$$

*with probability at least*
$$1 - \frac{2}{w}.$$

## 1.4 Additional Comments

As we have used different notions coming from different fields and different approaches, and we want to make our pages as readable as possible, we have included most notions and most basic facts in separate Sections. For instance, in Section 2 we introduce the basic notions and notations to deal with the Linear Algebra case, whereas in Section 4 we do the corresponding task for the Non–Linear case. Thus, we show Theorems 1 and 2 in Section 3 and Theorems 3 and 4 in Section 5. Finally, Section 6 is devoted to introduce the basic notions and notations to show Corollary 5 above.

# 2 Basic notions and notations

## 2.1 Notations for the real projective space

For every positive integer number $m \in \mathbb{N}$, let $\mathbb{P}_m(\mathbb{R})$ be the real projective space of dimension $m$. We denote by $\pi_\mathbb{R} : \mathbb{R}^{m+1} \setminus \{0\} \longrightarrow \mathbb{P}_m(\mathbb{R})$ the canonical projection onto the projective space. Namely, given $X := (x_0, \ldots, x_m) \in \mathbb{R}^{m+1} \setminus \{0\}$ we denote by $\pi_\mathbb{R}(X) \in \mathbb{P}_m(\mathbb{R})$ the projective point whose homogeneous coordinates are given by the following identity :

$$\pi_\mathbb{R}(X) := (x_0 : x_1 : \ldots : x_m) \in \mathbb{P}_m(\mathbb{R}).$$

Let $<\cdot,\cdot> \; : \mathbb{R}^{m+1} \times \mathbb{R}^{m+1} \longrightarrow \mathbb{R}$ be the canonical Euclidean inner product in $\mathbb{R}^{m+1}$. Namely, given $X := (x_0, \ldots, x_m), Y := (y_0, \ldots, y_m) \in \mathbb{R}^{m+1}$ we define

$$<X,Y> := \sum_{i=0}^m x_i y_i.$$

For every $X \in \mathbb{R}^{m+1}$ we denote its canonical Euclidean norm as $\|X\| := (<X,X>)^{1/2}$.

This Euclidean inner product induces a natural Riemannian structure on $\mathbb{P}_m(\mathbb{R})$ which we denote by $(\mathbb{P}_m(\mathbb{R}), can)$. The Riemannian metric on $(\mathbb{P}_m(\mathbb{R}), can)$ is denoted by $d_R$. The metric $d_R$ is given by the following rule : given $\pi_\mathbb{R}(X), \pi_\mathbb{R}(Y) \in \mathbb{P}_m(\mathbb{R})$, the Riemannian distance between $\pi_\mathbb{R}(X)$ and $\pi_\mathbb{R}(Y)$ is given by the following identity :

$$d_R(\pi_\mathbb{R}(X), \pi_\mathbb{R}(Y)) := \arccos \frac{\| <X,Y> \|}{\|X\|\|Y\|}.$$

We also introduce a Fubini–Study metric $d_{FS}$ on $\mathbb{P}_m(\mathbb{R})$ from the previous Riemannian metric in the following terms : given $\pi_\mathbb{R}(X), \pi_\mathbb{R}(Y) \in \mathbb{P}_m(\mathbb{R})$, the Fubini–Study distance from $\pi_\mathbb{R}(X)$ to $\pi_\mathbb{R}(Y)$ is given by the following identity :

$$d_{FS}(\pi_\mathbb{R}(X), \pi_\mathbb{R}(Y)) := \sin \, d_R(\pi_\mathbb{R}(X), \pi_\mathbb{R}(Y)).$$

For every $V \subseteq \mathbb{P}_m(\mathbb{R})$, let $\widetilde{V} \subset \mathbb{R}^{m+1}$ be the cone over $V$, namely

$$\widetilde{V} := \pi_\mathbb{R}^{-1}(V) \cup \{0\} \subseteq \mathbb{R}^{m+1}.$$

The Fubini–Study metric satisfies the following useful property.

**Lemma 6** *Let $V \subseteq \mathbb{P}_m(\mathbb{R})$ be a subset of the real projective space and let $\widetilde{V} := \pi_\mathbb{R}^{-1}(V) \cup \{0\} \subseteq \mathbb{R}^{m+1}$ be the corresponding cone over $V$. Then, for every point $X \in \mathbb{R}^{m+1}$ such that $\|X\| = 1$, the following equality holds :*

$$d(X, \widetilde{V}) := d_{FS}(\pi_\mathbb{R}(X), V) = inf\{d_{FS}(\pi_\mathbb{R}(X), \pi_\mathbb{R}(Y)) : \pi_\mathbb{R}(Y) \in V\},$$

*where $d(X, \widetilde{V}) := inf\{\|X - Z\| : Z \in \widetilde{V}\}$ and $\|\cdot\|$ is the canonical Euclidean norm.*

Let $S^m$ be the unit sphere in $\mathbb{R}^{m+1}$, namely

$$S^m := \{X \in \mathbb{R}^{m+1} : \|X\| = 1\}.$$

The canonical Euclidean inner product in $\mathbb{R}^{m+1}$ also induces a Riemannian structure on $S^m$ which we denote by $(S^m, can)$. Let us also denote by $p_\mathbb{R} : S^m \longrightarrow \mathbb{P}_m(\mathbb{R})$ the natural projection. Namely,

$$p_\mathbb{R} := \pi_\mathbb{R} \mid_{S^m}.$$

For every measurable subset $U \subseteq \mathbb{R}^{m+1}$ we denote by $Vol(U)$ the standard Lebesgue measure of $U$, i.e.

$$\int_{\mathbb{R}^{m+1}} \chi_U dx_1 \cdots dx_{m+1},$$

where $\chi_U : \mathbb{R}^{m+1} \longrightarrow \mathbb{R}$ is the characteristic function associated to $U$.
The canonical Riemannian structures on $\mathbb{P}_m(\mathbb{R})$ and $S^m$ respectively yield volume forms on every space. These volume forms lead to the following notions and notations.
For every measurable subset $U \subset S^m$ we define the spherical volume of $U$ in the following terms

$$Vol_S(U) := \int_{S^m} \chi_U(\theta) d\theta,$$

where $d\theta$ is the volume form in spherical coordinates. Observe that

$$Vol_S(S^m) = \frac{2\pi^{\frac{m+1}{2}}}{\Gamma(\frac{m+1}{2})},$$

where $\Gamma$ is the gamma function and $\chi_U : S^m \longrightarrow \mathcal{R}$ is the characteristic function associated to $U$.
Finally, for every measurable subset $U \subseteq \mathbb{P}_m(\mathbb{R})$, its projective volume satisfies :

$$Vol_P(U) := \frac{1}{2} Vol_S(p_\mathbb{R}^{-1}(U)).$$

## 2.2 Condition numbers for square matrices with real entries

As observed in the Introduction, the usual condition numbers of Linear Algebra are naturally defined in the projective space $\mathbb{P}(\mathcal{M}_n(\mathbb{R})) \cong \mathbb{P}_{n^2-1}(\mathbb{R})$. Here we recall these notions and some standard facts about them.
Let $S \subseteq \mathcal{M}_n(\mathbb{R})$ be the algebraic variety of all $n \times n$ singular matrices. Namely,

$$S := \{M \in \mathcal{M}_n(\mathbb{R}) : \det(M) = 0\}.$$

Let us identify $\mathcal{M}_n(\mathbb{R}) \cong \mathbb{R}^{n^2}$. The norm associated to the canonical Euclidean inner product in $\mathcal{M}_n(\mathbb{R})$ is usually called the Frobenius (also Hilbert–Weil, also Schur) norm. Given a square $n \times n$ real matrix $M \in \mathcal{M}_n(\mathbb{R})$ we denote its Frobenius norm as $\|M\|_F$. We also denote by $d_F(M, N)$ the Frobenius distance of two matrices $M, N \in \mathcal{M}_n(\mathbb{R})$. We denote by $\|M\|$ its norm as linear operator, namely

$$\|M\| := sup\{\|Mx\| : x \in \mathbb{R}^n \ \|x\| = 1\}.$$

An outstanding statements concerning the norm $\|M\|$ is the following statement due to C. Eckardt and G. Young.

**Theorem 7** *[6] For every $M \in \mathcal{M}_n(\mathbb{R})$,*

$$\|M^{-1}\| = \frac{1}{d_F(M, S)},$$

*where $d_F$ is the Frobenius distance.*

The usual condition number of a square matrix $M \in \mathcal{M}_n(\mathbb{R})$ is given by the following identity :
$$k(M) := \|M\|\|M^{-1}\|.$$

A second (and useful) condition number $\mu$ was discussed in [34, 3]. It can be defined as follows: for every $M \in \mathcal{M}_n(\mathbb{R})$, we define

$$\mu(M) := \|M^{-1}\|\|M\|_F.$$

As observed in [34, 3], the following inequality holds for every matrix $M \in \mathcal{M}_n(\mathbb{R})$ :
$$k(M) \leq \mu(M), \quad \forall M \in \mathcal{M}_n(\mathbb{R}). \tag{5}$$

For every matrix $M \in \mathcal{M}_n(\mathbb{R})$, both condition numbers $k(M)$ and $\mu(M)$ depend only on the projective point $\pi_{\mathbb{R}}(M) \in \mathbb{P}(\mathcal{M}_n(\mathbb{R}))$. Thus, we preserve the notation $k$ and $\mu$ to denote the mappings $k, \mu : \mathbb{P}(\mathcal{M}_n(\mathbb{R})) \longrightarrow \mathbb{R}$. Namely, $k(\pi_{\mathbb{R}}(M)) := k(M)$ and $\mu(\pi_{\mathbb{R}}(M)) := \mu(M)$. Let us denote by $\widehat{S} \subseteq \mathbb{P}(\mathcal{M}_n(\mathbb{R}))$ given by the following identity :

$$\widehat{S} := \pi_{\mathbb{R}}(S \setminus \{0\}).$$

For every projective point $\pi_{\mathbb{R}}(M) \in \mathbb{P}(\mathcal{M}_n(\mathbb{R}))$ we define the distance to the algebraic variety of singular matrices $\hat{S}$ as

$$\rho(\pi_{\mathbb{R}}(M)) := d_{FS}\left(\pi_{\mathbb{R}}(M), \widehat{S}\right) := inf\{d_{FS}(\pi_{\mathbb{R}}(M), Z) : Z \in \widehat{S}\}.$$

The following statement is an outstanding consequence of Eckardt & Young Theorem above.

**Corollary 8** For every $Z \in \mathbb{P}(\mathcal{M}_n(\mathbb{R}))$, the following holds

$$k(Z) \leq \mu(Z) = \frac{1}{\rho(Z)}.$$

*Proof.–* This fact is a consequence of Theorem 7 and Lemma 6 above. The ideas are essentially due to [34] (cf. also [3]). We include a proof for didactical reasons only.

The first inequality is Inequality 5 above. As for the equality, let $M \in \mathcal{M}_n(\mathbb{R})$ be a square matrix such that $\pi_R(M) = Z$. Let $M' \in S^{n^2-1} \subseteq \mathcal{M}_n(\mathbb{R})$ be the square matrix given by :

$$M' := \frac{1}{\|M\|_F} M \in \mathcal{M}_n(\mathbb{R}).$$

We obviously have

$$k(M') = k(M)$$
$$\mu(M') = \mu(M)$$

and $\pi_\mathbb{R}(M) = p_\mathbb{R}(M') = Z$.

Moreover, as $\|M'\|_F = 1$, we conclude from Theorem 7 above that the following equality holds :

$$\mu(M) = \mu(M') = \|M'^{-1}\| = \frac{1}{d_F(M', S)}.$$

Now, observe that $S$ is a cone. Namely, $S = \pi_\mathbb{R}^{-1}(\widehat{S}) \cup \{0\}$. Then, from Lemma 6 above, we conclude :

$$d_F(M', S) = d_{FS}\left(\pi_\mathbb{R}(M'), \widehat{S}\right) = \rho(Z),$$

and the statement follows. ∎

For every $\epsilon > 0$, let $\widehat{S}(\epsilon)$ be the compact neighborhood of $\widehat{S}$ in $\mathbb{P}(\mathcal{M}_n(\mathbb{R}))$ given in the following terms :

$$\widehat{S}(\epsilon) := \{Z \in \mathbb{P}(\mathcal{M}_n(\mathbb{R})) : d_{FS}(Z, \widehat{S}) \leq \epsilon\} = \{Z \in \mathbb{P}(\mathcal{M}_n(\mathbb{R})) : \rho(Z) \leq \epsilon\}.$$

The following statement is a mere consequence of Corollary 8 above.

**Corollary 9**

$$\widetilde{\widehat{S}(\varepsilon)} = \pi_\mathbb{R}^{-1}(\widehat{S}(\epsilon)) \cup \{0\} = \{M \in \mathcal{M}_n(\mathbb{R}) : \mu(M) \geq \frac{1}{\epsilon}\}.$$

The following Theorem can be found in [3].

**Theorem 10** *With the same notations and assumptions as above, the following holds for every $\varepsilon > 0$ :*

$$Vol_P(\widehat{S}(\varepsilon)) \leq \varepsilon Vol_P\left(\mathbb{P}(\mathcal{M}_n(\mathbb{R}))\right) n^{5/2}.$$

Sharper estimates can also be found in [7], [8] but they do not change essentially our forthcoming arguments. From this Theorem we may easily conclude the following statement.
Let $H, \varepsilon \in \mathbb{R}$ be two positive real numbers, $\varepsilon > 0$. Let $B(0, H) \subseteq \mathcal{M}_n(\mathbb{R})$ be the closed ball of radius $H$ centered at the origin with respect to the canonical Euclidean norm. We define the compact subset

$$R(\varepsilon, H) := B(0, H) \bigcap \widetilde{\widehat{S}(\varepsilon)} = B(0, H) \bigcap \left(\pi_\mathbb{R}^{-1}(\widehat{S}(\varepsilon)) \cup \{0\}\right). \qquad (6)$$

**Corollary 11** *Let $H, \varepsilon \in \mathbb{R}$ be two positive real numbers, $\varepsilon > 0$. With the previous notations and assumptions, the following inequality holds :*

$$Vol\left(R(\varepsilon, H)\right) \leq \varepsilon n^{1/2} Vol_S(S^{n^2-1}) H^{n^2}.$$

*Proof.–* Let us denote by $R$ the compact set $R(\varepsilon, H)$. By integration in spherical coordinates we obtain the following equality.

$$Vol(R) := \int_{\mathbb{R}^{n^2}} \chi_R dx_{1,1} \cdots dx_{n,n} = \int_{S^{n^2-1}} \int_0^H \chi_R(r\theta) r^{n^2-1} dr d\theta.$$

Now, let us observe that $\chi_R$ is homogeneous of degree 0 (for $0 \leq r \leq H$). Namely, given $r, r' \in \mathbb{R}$, $0 \leq r, r' \leq H$, $\chi_R(r\theta) = \chi_R(r'\theta)$, for all $\theta$.
Then, since $\pi_R(R) = \widehat{S}(\varepsilon)$ we conclude the following identity :

$$Vol(R) := \left(\int_{S^{n^2-1}} \chi_R(\theta) d\theta\right) \left(\int_0^H r^{n^2-1} dr\right) = \frac{H^{n^2}}{n^2} Vol_S(p_\mathbb{R}^{-1}(\widehat{S}(\varepsilon))),$$

and the statement follows from Theorem 10 above. ■

## 2.3 Davenport's discrepancy Bounds

Here we recall a beautiful statement due to H. Davenport (cf. [5]). Let $m \in \mathbb{N}$ be a positive integer number.
Given a measurable subset $R \subseteq \mathbb{R}^m$ of finite volume $Vol(R) < \infty$, we denote by $N(R)$ the number of points in $R$ with integer coordinates. Namely,

$$N(R) := \sharp\left(R \bigcap \mathbb{Z}^m\right).$$

It is well–known (cf. [20], for instance) that the following equality holds :

$$\lim_{h \to \infty} \frac{N(hR)}{h^m} = Vol(R),$$

where $h \in \mathbb{R}$ is a positive real number and $hR := \{hx \; : \; x \in R\}$. Upper bounds for $\|N(hR) - Vol(hR)\|$ are usually called *discrepancy bounds* (cf. [10], [11], [37] and references therein for other approaches). Here we recall a method due to H. Davenport to obtain discrepancy bounds which by the way is very accurate for the problems we deal with.

For every positive integer number $\ell \in \mathbb{N}$, $1 \leq \ell \leq m$, we denote by $\mathcal{I}_\ell^{(m)}$ the class of all subsets $I$ of $\{i \; : \; 1 \leq i \leq m\}$ such that $\sharp(I) = \ell$.

For every subset $I \in \mathcal{I}_\ell^{(m)}$, let us denote by $\Phi_I$ the projection

$$\Phi_I \; : \mathbb{R}^m \longrightarrow \mathbb{R}^\ell,$$

given by

$$\Phi_I(x_1, \ldots, x_m) := (x_i \; : \; i \in I).$$

Let $R$ be a compact subset of $\mathbb{R}^m$. We denote by $V(R, \ell)$ the following quantity :

$$V(R, \ell) := \sum_{I \in \mathcal{I}_\ell^{(m)}} Vol_\ell(\Phi_I(R)),$$

where $Vol_\ell(\Phi_I(R))$ denotes the Lebesgue of $\Phi_I(R)$ as subset of $\mathbb{R}^\ell$. We define $V(R, 0) = 1$ by convention.

For every subset $U \subseteq \mathbb{R}^m$, we denote by $\beta_0(U)$ the number of connected components of $U$ (i.e. the 0−th Betti number).

Let $I \in \mathcal{I}_\ell^{(m)}$ be given, and let $J \subseteq I$, such that $\sharp(I \setminus J) = 1$. For every point $A := (\alpha_i \; : \; i \in J) \in \mathbb{R}^{\ell-1}$, let $r_J(A)$ be the real line in $\mathbb{R}^\ell$ given by the following identity :

$$r_J(A) := \{(x_i : i \in I) : x_k = \alpha_k, \forall k \in J\}.$$

Let $R \subseteq \mathbb{R}^{n^2}$ be a compact set. We define $h(R, \ell)$ as the following number :

$$h(R, \ell) := sup\{\beta_0(\Phi_I(R) \cap r_J(A)) : I \in \mathcal{I}_\ell^{(m)}, J \subseteq I, \sharp(I \setminus J) = 1, A \in \mathbb{R}^{\ell-1}\},$$

where $\beta_0$ denotes the 0-th Betti number (i.e. the number of connected components.

We finally define $h(R) := sup\{h(R, \ell) : 1 \leq \ell \leq m - 1\}$.

These notations stated, the following statement holds.

**Theorem 12 (Davenport's discrepancy bound)** *[5] Let $R \subseteq \mathbb{R}^m$ a compact set. Let $N(R)$ be the number of points in $R$ with integer coordinates. Then, the following holds :*

$$|N(R) - Vol(R)| \leq \sum_{\ell=0}^{m-1} h(R)^{m-\ell} V(R, \ell).$$

## 2.4 Northcott–Schmidt's height, bit length and visible points from the origin

Let $\Lambda \subseteq \mathbb{R}^{m+1}$ be a lattice. A non–zero point $X \in \Lambda \setminus \{0\}$ is said to be *visible from the origin* (or simply *visible*) if there is no point of $\Lambda$ in the segment $[0, X] \subseteq \mathbb{R}^{m+1}$ between the origin and $X$.

For every non–zero real number $\rho \in \mathbb{R} \setminus \{0\}$, let $\mathbb{Z}^{m+1}(\rho)$ be the lattice given by the following identity :

$$\mathbb{Z}^{m+1}(\rho) := \{\rho X \ : \ X \in \mathbb{Z}^{m+1}\}.$$

It is easy to prove that the set of visible points in $\mathbb{Z}^{m+1}(\rho)$ is the set of all those points $\rho X := (\rho x_0, \ldots, \rho x_m)$ such that the greatest common divisor $gcd(x_0, \ldots, x_m) = 1$.

For instance, let $\mathcal{M}_n(\mathbb{Z})$ be the set of all $n \times n$ matrices with integer entries. The set $\mathcal{M}_n(\mathbb{Z})$ is a lattice in $\mathcal{M}_n(\mathbb{R}) \cong \mathbb{R}^{n^2}$. For every non–zero real number $\rho \in \mathbb{R}$, let $\mathcal{M}_n(\mathbb{Z})(\rho)$ be the lattice in $\mathcal{M}_n(\mathbb{R})$ given by the following identity :

$$\mathcal{M}_n(\mathbb{Z})(\rho) := \{\rho M \ : \ M \in \mathcal{M}_n(\mathbb{Z})\}.$$

The visible points of $\mathcal{M}_n(\mathbb{Z})(\rho)$ are exactly those matrices $\rho M := (\rho x_{i,j})_{1 \leq i,j \leq n}$ such that the greatest common divisor $gcd(x_{i,j} \ : \ 1 \leq i, j \leq n) = 1$.

Visible points in a lattice $\mathbb{Z}^{m+1}(\rho)$ and projective points are closely related. In order to illustrate this relation, let us introduce the *Schmidt height* of a projective point. Let $P \subseteq \mathbb{N}$ be the class of all prime numbers. For every $p \in P$, let $\|\cdot\|_p \ : \ \mathbb{Q} \longrightarrow \mathbb{R}_+$ be the non–archimedean $p$–adic absolute value.

**Definition 13** *Let $\pi_\mathbb{R}(X) := (x_0 : \ldots : x_m) \in \mathbb{P}_m(\mathbb{Q})$ a projective point whose homogeneous coordinates are rationals, i.e. assume that $X = (x_0, \ldots, x_m) \in \mathbb{Q}^{m+1} \setminus \{0\}$. We define the Northcott–Schmidt absolute height of $\pi_\mathbb{R}(X)$ as the following infinite product :*

$$H(\pi_\mathbb{R}(X)) := \|X\| \prod_{p \in P} \max\{\|x_i\|_p : 0 \leq i \leq m\}.$$

This notion is well–defined because of *Weil's product formula*. Let us consider now the projective space $\mathbb{P}(\mathcal{M}_n(\mathbb{Q}))$ define by the vector space of all $n \times n$ matrices with rational entries $\mathcal{M}_n(\mathbb{Q})$. Given a non–zero matrix $M := (x_{i,j}) \in \mathcal{M}_n(\mathbb{Q})$, the Northcott–Schmidt absolute height of $\pi_\mathbb{R}(M)$ is given by the following identity :

$$H(\pi_\mathbb{R}(M)) := \|M\|_F \prod_{p \in P} \max\{\|x_{i,j}\|_p : 1 \leq i, j \leq n\},$$

where $\|M\|_F$ is the Frobenius norm introduced above.

Let $\rho X := (\rho x_0, \ldots, \rho x_m) \in \mathbb{Z}^{m+1}(\rho)$ be a visible point in the lattice $\mathbb{Z}^{m+1}(\rho)$. Let $\pi_\mathbb{R}(\rho X) = \pi_\mathbb{R}(X) \in \mathbb{P}_m(\mathbb{Q})$ the corresponding projective point. Then, it is easy to show the following Lemma.

**Lemma 14** *With the same notations and assumptions as above, the following properties hold :*

i) *Let $\rho X := (\rho x_0, \ldots, \rho x_m) \in \mathbb{Z}^{m+1}(\rho)$ be a point in the lattice $\mathbb{Z}^{m+1}(\rho)$. Then, $\rho X$ is a visible point from the origin if and only if*

$$\|X\| = min\{\|Y\| : \rho Y \in \mathbb{Z}^{m+1}(\rho), \pi_\mathbb{R}(\rho Y) = \pi_\mathbb{R}(X)\}.$$

ii) *Let $\rho X \in \mathbb{Z}^{m+1}(\rho)$ be a visible point. Then,*

$$\|X\| = H(\pi_\mathbb{R}(\rho X)).$$

iii) *Let $V \subseteq \mathbb{P}(\mathcal{M}_n(\mathbb{R}))$ be a subset of a real projective space and let $\widetilde{V} := \pi^{-1}(V) \cup \{0\} \subseteq \mathcal{M}_n(\mathbb{R})$ the cone over $V$. Let $B(0, H)$ be the closed ball in $\mathbb{R}^{m+1}$ of radius $H$ centered at the origin. Then, the following holds :*

- *For every $\pi_\mathbb{R}(X) \in V \cap \mathbb{P}_m(\mathbb{Q})$ there is a visible point $Y \in \widetilde{V} \cap \mathbb{Z}^{m+1}$ such that $\pi_\mathbb{R}(Y) = \pi_\mathbb{R}(X)$ and $\|Y\| = H(\pi_\mathbb{R}(X))$.*
- *The number of visible points of $\mathbb{Z}^{m+1}$ that lie in $\widetilde{V} \cap B(0, H)$ is two times the number of projective points in $V \cap \mathbb{P}_m(\mathbb{Q})$ of Northcott-Schmidt absolute height at most $H$.*

In fact, given a projective point $\pi_\mathbb{R}(X) \in \mathbb{P}_m(\mathbb{Q})$ and given a visible point $Y \in \mathbb{Z}^{m+1}$ such that $\pi_\mathbb{R}(Y) = \pi_\mathbb{R}(X)$, the two visible points in $\mathbb{Z}^{m+1}$ associated to $\pi_\mathbb{R}(X)$ are $\{Y, -Y\} \subseteq \mathbb{Z}^{m+1}$.

Given a projective point $\pi_\mathbb{R}(X) \in \mathbb{P}_m(\mathbb{Q})$, we may represent this point by some visible point $Y \in \mathbb{Z}^{m+1}$ such that $\pi_\mathbb{R}(Y) = \pi_\mathbb{R}(X)$. Thus the bit length of $\pi_\mathbb{R}(X)$ is the number of tapes cells required to encode $Y$ in a Turing machine. As $Y \in \mathbb{Z}^{m+1}$, this quantity is essentially equal to the logarithm of the norm of $Y$, i.e. $\log \|Y\|$. This justifies the following Definition.

**Definition 15** *Given a projective point $\pi_\mathbb{R}(X) \in \mathbb{P}_m(\mathbb{Q})$, we define its bit length as the logarithm of its Northcott-Schmidt absolute height. Namely,*

$$bl(\pi_\mathbb{R}(X)) = \log_2 H(\pi_\mathbb{R}(X)).$$

# 3 The Linear Algebra case

## 3.1 Davenport's discrepancy in the Linear Algebra case

In this Subsection we apply Davenport's strategy (as introduced in Subsection 2.3 above) to count the number of integer points in the region $R(\varepsilon, H)$ introduced in Equation (6) of Subsection 2.2 above. Namely, given two positive real numbers $H, \varepsilon \in \mathbb{R}$, let $R(\varepsilon, H) \subseteq \mathcal{M}_n(\mathbb{R})$ be the compact subset given by the following identity :

$$R(\varepsilon, H) := B(0, H) \bigcap \widetilde{\widehat{S}(\varepsilon)} = B(0, H) \bigcap \left( \pi_\mathbb{R}^{-1}(\widehat{S}(\varepsilon)) \cup \{0\} \right).$$

**Proposition 16** *With the previous notations and assumptions, let $H \geq 1$ be a real number and let $N(\varepsilon, H)$ be the number of matrices with integer entries in $R(\varepsilon, H)$, namely*

$$N(\varepsilon, H) := \sharp \left( R(\varepsilon, H) \bigcap \mathcal{M}_n(\mathbb{Z}) \right).$$

*Then, the following inequality holds :*

$$|N(\varepsilon, H) - Vol\left(R(\varepsilon, H)\right)| \leq T_n \mathfrak{S}^{(n^2)} H^{n^2-1}.$$

*where $\mathfrak{S}^{(n^2)}$ is the constants introduced in Equation 1 of the Introduction and*

$$T_n := (2 \max\{4, n\})^{(2n^4+4n^2)},$$

*is the constant introduced in Theorem 1.*
*In particular, from Corollary 11 we conclude the following inequality*

$$N(\varepsilon, H) \leq \varepsilon n^{1/2} Vol_S(S^{n^2-1}) H^{n^2} + T_n \mathfrak{S}^{(n^2)} H^{n^2-1}.$$

*Proof.–* First of all, observe that $R(\varepsilon, H)$ is the set of all matrices $M \in \mathcal{M}_n(\mathbb{R})$ that satisfy the following properties :

- $\|M\|_F^2 \leq H^2$,
- $d_{FS}(\pi_\mathbb{R}(M), \widehat{S}) \leq \varepsilon$.

Now, let $M' \in \mathcal{M}_n(\mathbb{R})$ be the square matrix given by the following identity :

$$M' := \frac{1}{\|M\|_F} M.$$

We obviously have, $\|M'\|_F = 1$ and $\pi_\mathbb{R}(M') = \pi_R(M)$. From Proposition 6 above, we conclude :

$$d_{FS}(\pi_\mathbb{R}(M), \widehat{S}) = d_{FS}(\pi_\mathbb{R}(M'), \widehat{S}) = d_F(M', S).$$

Then, the previous two properties that define $R(\varepsilon, H)$ may be replaced by the next two ones : For every $M \in \mathcal{M}_n(\mathbb{R})$, $M \in R(\varepsilon, H)$ if and only the following two properties hold :

- $\|M\|_F^2 \leq H^2$,

- $\exists N \in S$, such that $\|M - \|M\|_F N\|_F \leq \varepsilon \|M\|_F$.

Then, with the notations of Subsection 2.3 above, for every integer number $\ell$, $1 \leq \ell \leq n^2$, let $\mathcal{I}_\ell^{(n^2)}$ be the class of all subsets $I$ of $\{(i,j) : 1 \leq i,j \leq n\}$ such that $\sharp(I) = \ell$. For every $I \in \mathcal{I}_\ell^{(n^2)}$, let $\Phi_I : \mathcal{M}_n(\mathbb{R}) \longrightarrow \mathbb{R}^\ell$ the projection defined by $I$ in the following terms :

$$\Phi_I\left((x_{i,j})_{1 \leq i,j \leq n}\right) := (x_{r,s} : (r,s) \in I) \in \mathbb{R}^\ell.$$

Observe that for every $I \in \mathcal{I}_\ell^{(n^2)}$ the following holds :

$$\Phi_I(R(\varepsilon, H)) \subset B_\ell(0, H),$$

where $B_\ell(0, H)$ is the closed ball of radius $H$ centered at the origin of $\mathbb{R}^\ell$. Then if $Vol(\Phi_I(R(\varepsilon, H)))$ denotes the $\ell-$volume of $\Phi_I(R(\varepsilon, H))$, the following inequality holds :

$$Vol(\Phi_I(R(\varepsilon, H))) \leq Vol(B_\ell(0, H)) = H^\ell K_\ell.$$

With the same notations as in Subsection 2.3 and Theorem 12 above, the following inequality holds :

$$V(R(\varepsilon, H), \ell) = \sum_{I \in \mathcal{I}_\ell^{(n^2)}} Vol(\Phi_I(R(\varepsilon, H))) \leq H^\ell \binom{n^2}{\ell} K_\ell.$$

The second invariant in Davenport's Theorem 12 is given in the following terms. Let $\ell \in \mathbb{N}$ be a positive integer number such that $1 \leq \ell \leq n^2$. Let $I \in \mathcal{I}_\ell^{(n^2)}$ be a subset of $\{(i,j) : 1 \leq i,j \leq n\}$ of cardinal $\ell$. Let $J \subseteq I$ a subset such that $\sharp(I \setminus J) = 1$. Let $A = (\alpha_{(i,j)})_{(i,j) \in J} \in \mathbb{R}^{\ell-1}$ be any affine point and let $r_J(A)$ be the real line in $\mathbb{R}^\ell$ given by

$$r_J(A) := \{(x_{i,j} : (i,j) \in I) \in \mathbb{R}^\ell : x_{i,j} = \alpha_{i,j}, \forall (i,j) \in J\}.$$

We define $h(R(\varepsilon, H), \ell)$ as the maximum number of connected components of the sets given by

$$\Phi_I(R(\varepsilon, H)) \cap r_J(A),$$

where $I \subseteq \mathcal{I}_\ell^{(n^2)}$, $J \subseteq I$ such that $\sharp(I \setminus J) = 1$ and $A \in \mathbb{R}^{\ell-1}$.

With the same notations, let us introduce the real algebraic subset $V(\varepsilon, H, I, J, A) \subset \mathbb{R}^{2n^2+4}$ given by the following property : a point $(M, N, t_1, t_2, t_3, t_4) \in \mathcal{M}_n(\mathbb{R})^2 \times \mathbb{R}^4 \cong \mathbb{R}^{2n^2+4}$ belongs to $V(\varepsilon, H, I, J, A)$ if and only if the following properties hold :

i) $x_{i,j} = \alpha_{ij}$, for all $(i,j) \in J$,

ii) $H^2 - \|M\|_F^2 = t_1^2$,

iii) $\det(N) = 0$,

iv) $\|M\|_F^2 = t_2^2$,

v) $t_2 - t_3^2 = 0$,

vi) $\|M\|_F^2 \varepsilon^2 - \|M - t_2 N\|_F^2 = t_4^2$

We observe that the maximum of the degrees involved in these equations is $\max\{4, n\}$. In other words, the polynomials of maximal degree occur either in the equation *iii)* or in equation *vi)*. Now we may apply the Milnor–Thom–Oleinik–Petrovsky upper bounds for the sum of the Betti numbers of real algebraic varieties (cf. [19], [35], [23], [24]). These upper bounds yield the following inequality :

$$\sum_{i \in \mathbb{N}} \beta_i \left( V(\varepsilon, H, I, J, A) \right) \leq (2 \max\{4, n\})^{2n^2 + 4}. \tag{7}$$

On the other hand, let $\Psi : \mathcal{M}_n(\mathbb{R})^2 \times \mathbb{R}^4 \longrightarrow \mathcal{M}_n(\mathbb{R})$ be the projection given by $\Psi(M, N, t_1, t_2, t_3, t_4) = M$. Then, the following identity holds :

$$\Phi_I \left( \Psi \left( V(\varepsilon, H, I, J, A) \right) \right) = \Phi_I(R(\varepsilon, H)) \cap r_J(A). \tag{8}$$

In order to see this identity, let us observe that Equation *ii)* above means $\|M\|_F \leq H$, Equation *iii)* means that $N$ is a singular matrix and, Equations *iv), v)* and *vi)* mean that $\|M - \|M\|_F N\| \leq \|M\|_F \varepsilon$.

Finally, as $\Psi$ and $\Phi_I$ are continuous mappings, we conclude from Identity (8) and Inequality (7) above the following inequality :

$$\beta_0 \left( \Phi_I(R(\varepsilon, H)) \cap r_J(A) \right) \leq \beta_0(V(\varepsilon, H, I, J, A)) \leq (2 \max\{4, n\})^{2n^2 + 4}.$$

In particular, for every $\ell \in \mathbb{N}$, $1 \leq \ell \leq n^2$, the following inequality holds :

$$h(R(\varepsilon, H), \ell) \leq (2 \max\{3, n\})^{2n^2 + 4}.$$

Finally, as $H \geq 1$ we conclude from Theorem 12 the following inequality :

$$|N(R(\varepsilon, H)) - Vol(R(\varepsilon, H))| \leq H^{n^2 - 1} (2 \max\{3, n\})^{(2n^4 + 4n^2)} \mathfrak{S}^{(n^2)}.$$

∎

**Proposition 17** *With the same notations and assumptions as in Proposition 16 above, let $H \geq 1$ be a positive real number. Let $N(1, H)$ be the number of matrices with integer entries in $R(1, H)$. Namely,*

$$N(1, H) := \sharp \left( B(0, H) \bigcap \mathcal{M}_n(\mathbb{Z}) \right).$$

Then, the following holds :

$$|N(1,H) - Vol\left(B(0,H)\right)| \leq H^{n^2-1}\mathfrak{S}^{(n^2)}.$$

In particular, the following inequality also holds :

$$\frac{Vol_S(S^{n^2-1})}{n^2}H^{n^2} - H^{n^2-1}\mathfrak{S}^{(n^2)} \leq N(1,H).$$

*Proof.–* The proof follows the same arguments as in Proposition 16 above, noting that $h(B(0,H),\ell)$ is always 1. ∎

## 3.2 On ill–conditioned matrices that are visible from the origin.

**Proposition 18** *Let $H, \varepsilon \in \mathbb{R}$ be two positive real numbers, $H \geq 1$. Let $\mathcal{N}(\varepsilon, H)$ be the number of points in $\mathbb{P}(\mathcal{M}_n(\mathbb{Q}))$ of absolute Northcott–Schmidt height at most $H$ which belong to $S(\varepsilon)$. Then, the following holds*

$$|\mathcal{N}(\varepsilon, H) - \frac{Vol(R(\varepsilon, H))}{2\zeta(n^2)}| \leq L(\varepsilon, n)H^{n^2-1} + \frac{H}{2},$$

*where*

$$L(\varepsilon, n) := \left(\frac{T_n\mathfrak{S}^{(n^2)}}{2\zeta(n^2-1)} + \frac{Vol(R(\varepsilon, 1))}{2}\right).$$

*and $T_n$ and $K_\ell^{(n^2)}$ are as in Proposition 16 above.*
*In particular, the following inequality is a consequence of Corollary 11 :*

$$\mathcal{N}(\varepsilon, H) \leq \frac{\varepsilon n^{1/2} Vol_S(S^{n^2-1})}{2\zeta(n^2)}H^{n^2} + L'(\varepsilon, n)H^{n^2-1} + \frac{H}{2},$$

*where*

$$L'(\varepsilon, n) := \left(\frac{T_n\mathfrak{S}^{(n^2)}}{2\zeta(n^2-1)} + \frac{\varepsilon n^{5/2} K_{n^2}}{2}\right).$$

*Proof.–* This proof is widely inspired by the Proof of Theorem 459 of [14]. Assume that $\varepsilon > 0$ is fixed throughout this Proof. Let $V(\varepsilon) \subseteq \mathcal{M}_n(\mathbb{R}) \cong \mathbb{R}^{n^2}$ be the compact subset given by the following identity :

$$V(\varepsilon) := B(0,1) \bigcap \left(\pi^{-1}\left(S(\varepsilon)\right) \cup \{0\}\right).$$

From Corollary 9, we easily conclude that for every positive real number $t \in \mathbb{R}_+$

$$tV(\varepsilon) := \{tX \ : \ X \in V(\varepsilon)\} = R(\varepsilon, t),$$

where $R(\varepsilon, t)$ is the compact set introduced in Subsection 3.1 above. In particular, for every positive real number $t \in \mathbb{R}_+$ the following identity holds :
$$t^{n^2} Vol(V(\varepsilon)) := Vol(R(\varepsilon, t)) \qquad (9)$$

For every positive real number $\rho \in \mathbb{R}_+$, let us denote by $A_\rho$ the finite subset of $B(0,1)$ given by the following equality
$$A_\rho := V(\varepsilon) \bigcap \mathcal{M}_n(\mathbb{Z})(\rho) \setminus \{0\},$$

where $\mathcal{M}_n(\mathbb{Z})(\rho)$ is the lattice in $\mathcal{M}_n(\mathbb{R})$ introduced in Subsection 2.4. Namely,
$$\mathcal{M}_n(\mathbb{Z})(\rho) := \{\rho M \ : \ M \in \mathcal{M}_n(\mathbb{Z})\}.$$

Let $g(\rho)$ be the number of points in $A_\rho$. Observe that the following equality holds :
$$g(\rho) := \sharp(A_\rho) = \sharp\left(R(\varepsilon, \rho^{-1}) \bigcap \mathcal{M}_n(\mathbb{Z})\right) - 1 = N(\varepsilon, \rho^{-1}) - 1. \qquad (10)$$

Let $f(\rho)$ be the number of points in $A_\rho$ which are visible from the origin. As in the Proof of Theorem 459 of [14] we may easily conclude that
$$g(\rho) := \sum_{i=1}^{\infty} f(m\rho).$$

From Möbius inversion formula (see Theorem 270 in [14], for instance) it follows that
$$f(\rho) := \sum_{m=1}^{\infty} \mu(m) g(m\rho),$$

where $\mu$ is Möbius function. Observe that Riemann's zeta function satisfies the following identity for every $s > 1$ (cf. Theorem 287 of [14], for instance) :
$$\frac{1}{\zeta(s)} := \sum_{m=1}^{\infty} \frac{\mu(m)}{m^s}. \qquad (11)$$

So
$$\rho^{n^2} f(\rho) - \frac{Vol(V(\varepsilon))}{\zeta(n^2)} = \sum_{m=1}^{\infty} \frac{\mu(m)}{m^{n^2}} \left((m\rho)^{n^2} g(m\rho) - Vol(V(\varepsilon))\right). \qquad (12)$$

Let $\mathbb{N}_\rho$ be the set of integer numbers given by the following identity
$$\mathbb{N}_\rho := \{m \in \mathbb{N} \ : \ 1 \leq m \leq \rho^{-1}\}.$$

Let us denote by $S_\rho$ the sum of the terms of the left hand series of Equation (30) whose indices are in $\mathbb{N}_\rho$. Namely,
$$S_\rho := \sum_{m \in \mathbb{N}_\rho} \frac{\mu(m)}{m^{n^2}} \left((m\rho)^{n^2} g(m\rho) - Vol(V(\varepsilon))\right).$$

Now, assume that $m \in \mathbb{N}_\rho$ (i.e. $m\rho \leq 1$). Then, replacing $(m\rho)^{-1}$ by $H$, we conclude from Identities (9) and (10) above the following identity :

$$|(m\rho)^{n^2} g(m\rho) - Vol(V(\varepsilon))| = \frac{1}{H^{n^2}}|N(\varepsilon, H) - (1 + Vol(R(\varepsilon, H)))| \quad (13)$$

From Proposition 16 above, we have :

$$\frac{1}{H^{n^2}}|N(\varepsilon, H) - Vol(R(\varepsilon, H))| \leq \frac{T_n \mathfrak{S}^{(n^2)}}{H}.$$

Replacing back $H$ by $(m\rho)^{-1}$, Equation (13) becomes the following inequality :

$$|(m\rho)^{n^2} g(m\rho) - Vol(V(\varepsilon))| \leq T_n \mathfrak{S}^{(n^2)} m\rho + (m\rho)^{n^2} \quad (14)$$

We conclude :

$$S_\rho \leq \sum_{m \in \mathbb{N}_\rho} \left( \frac{\mu(m)}{m^{n^2-1}} \rho T_n \mathfrak{S}^{(n^2)} + \mu(m)\rho^{n^2} \right). \quad (15)$$

Then,

$$S_\rho \leq \frac{\rho T_n \mathfrak{S}^{(n^2)}}{\zeta(n^2-1)} + \left( \sum_{m \in \mathbb{N}_\rho} \mu(m)\rho^{n^2} \right) \leq \frac{\rho T_n \mathfrak{S}^{(n^2)}}{\zeta(n^2-1)} + \rho^{n^2-1}. \quad (16)$$

On the other hand, let $\mathbb{N}'_\rho$ be the set of positive integer numbers which are not in $\mathbb{N}_\rho$. Namely,

$$\mathbb{N}'_\rho := \{m \in \mathbb{N} \ : \ m\rho > 1\}.$$

Now, assume $m \in \mathbb{N}'_\rho$ and $X \in A_{m\rho}$. Then, we have $X = (m\rho)Y$, where $Y \in \mathbb{Z}^{n^2}$. In particular, the following inequality holds :

$$\|Y\| = (m\rho)\|Y\| \geq m\rho > 1.$$

In other words, for every $m \in \mathbb{N}'_\rho$, $g(m\rho) = 0$. Hence, let us define $S'_\rho$ as the sum of the terms in the left hand series of Equation (12) above. Namely :

$$S'_\rho := \sum_{m \in \mathbb{N}'_\rho} \frac{\mu(m)}{m^{n^2}} |(m\rho)^{n^2} g(m\rho) - Vol(V(\varepsilon))|.$$

We have :

$$S'_\rho \leq \sum_{m \in \mathbb{N}'_\rho} \frac{\mu(m)}{m^{n^2}} |Vol(V(\varepsilon))| \leq V(\varepsilon) \left( \sum_{m \geq \rho^{-1}} \frac{1}{m^{n^2}} \right) \leq \rho Vol(V(\varepsilon)) \quad (17)$$

Now we combine Inequalities (16) and (17) to conclude :

$$|\rho^{n^2} f(\rho) - \frac{Vol(V(\varepsilon))}{\zeta(n^2)}| \leq S_\rho + S'_\rho \leq \frac{\rho T_n \mathfrak{S}^{(n^2)}}{\zeta(n^2-1)} + \rho^{n^2-1} + \rho Vol(V(\varepsilon)).$$

Finally, replacing $\rho^{-1}$ by $H$ in this Equation and using the identity described in Equation (9) we conclude :

$$|f(H^{-1}) - \frac{Vol(R(\varepsilon, H))}{\zeta(n^2)}| \leq \frac{\mathfrak{S}^{(n^2)}}{\zeta(n^2-1)} H^{n^2-1} + H + \frac{Vol(R(\varepsilon, H))}{H}.$$

Noting that $f(H^{-1}) = 1/2\mathcal{N}(\varepsilon, H)$ (cf. Lemma 14 above) the Proposition is achieved. ∎

**Proposition 19** *Let $H \in \mathbb{R}$ be a positive real number, $H \geq 1$. Let $\mathcal{N}(1, H)$ be the number of points in $\mathbb{P}(\mathcal{M}_n(\mathbb{Q}))$ of Northcott–Schmidt absolute height at most $H$. Then, the following holds*

$$|\mathcal{N}(1, H) - \frac{Vol(B(0, H))}{2\zeta(n^2)}| \leq L(1, n) H^{n^2-1} + \frac{H}{2},$$

*where*

$$L(1, n) := \frac{\mathfrak{S}^{(n^2)}}{2\zeta(n^2-1)} + \frac{K_{n^2}}{2}.$$

*where $T_n$ and $K_\ell^{(n^2)}$ are as in Proposition 16 above, and $\zeta$ is Riemann's zeta function. In particular, the following inequality also holds :*

$$\frac{Vol_S(S^{n^2-1})}{2n^2 \zeta(n^2)} H^{n^2} - \left( L(1, n) H^{n^2-1} + \frac{H}{2} \right) \leq \mathcal{N}(1, H).$$

*Proof.–* This is step by step the same proof as that of Proposition 18 using Proposition 17 instead of Proposition 16. ∎

### 3.3 Proof of Theorem 1

*Proof of Theorem 1.–* First of all, observe that the probability that a random choice of a system of linear equations $M \in \mathbb{P}(\mathcal{M}_n(\mathbb{Q}))$ of Nosthcott–Schmidt height at most $H$ satisfies $\mu(M) \geq 1/\varepsilon$ is given by the following quotient :

$$\frac{\mathcal{N}(\varepsilon, H)}{\mathcal{N}(1, H)}.$$

Let $A$ be the quantity :

$$A := \frac{K_{n^2}}{2\zeta(n^2)}.$$

From Propositions 18 and 19 we conclude :

$$\frac{\mathcal{N}(\varepsilon, H)}{\mathcal{N}(1, H)} \leq \frac{\varepsilon n^{5/2} A H + L'(\varepsilon, n) + \frac{1}{2H^{n^2-2}}}{AH - \left[ L(1, n) + \frac{1}{2H^{n^2-2}} \right]}.$$

Hence, we have

$$\frac{\mathcal{N}(\varepsilon, H)}{\mathcal{N}(1, H)} \leq \frac{\varepsilon n^{5/2} H + \frac{L'(\varepsilon, n)}{A} + \frac{1}{2AH^{n^2-2}}}{H - \left[\frac{L(1,n)}{A} + \frac{1}{2AH^{n^2-2}}\right]}.$$

Finally, as $H \geq 1$ we conclude :

$$\frac{\mathcal{N}(\varepsilon, H)}{\mathcal{N}(1, H)} \leq \varepsilon n^{5/2} + \frac{\varepsilon n^{5/2} \left[\frac{L(1,n)}{A} + \frac{1}{2AH^{n^2-2}}\right] + \frac{L'(\varepsilon, n)}{A} + \frac{1}{2AH^{n^2-2}}}{H - \left[\frac{L(1,n)}{A} + \frac{1}{2AH^{n^2-2}}\right]}.$$

where

$$B(\varepsilon, n) := \varepsilon n^{5/2} \left[\frac{L(1,n)}{A} + \frac{1}{2A}\right] + \frac{L'(\varepsilon, n)}{A} + \frac{1}{2A},$$

and

$$C(n) := \left[\frac{L(1,n)}{A} + \frac{1}{2A}\right].$$

∎

## 4 Condition Numbers for Systems of Multivariate Polynomial Equations

### 4.1 Two Hermitian inner products in $\mathcal{H}_{(d)}$.

For the sake of readability, we recall here some of the notations stated in Subsection 1.2 of the Introduction. We define $H_d$ as the set of all complex homogeneous polynomials in $n+1$ variables of degree $d$. Namely,

$$H_d := \{f \in \mathbb{C}[X_0, \ldots, X_n] \ : \ f \text{ homogeneous}, \deg(f) = d\}.$$

We define the standard Hermitian inner product in $H_d$ by identifying $H_d \cong \mathbb{C}^{N_d}$, where $N_d$ is the number of coefficients of a generic homogeneous polynomial $f \in \mathbb{C}[X_0, \ldots, X_N]$ of degree $d$. We have

$$N_d := \binom{d+n}{n}.$$

The standard Hermitian inner product $<\cdot, \cdot> \ : H_d \times H_d \longrightarrow \mathbb{C}$ is given in the following terms : Given $f, g \in H_d$ we define $<f, g>$ as

$$<f, g> := \sum_{|\mu|=d} a_\mu \overline{b_\mu},$$

where $\mu := (\mu_0, \ldots, \mu_n) \in \mathbb{N}^{n+1}$ is a multi–index, $|\mu| := \mu_0 + \cdots + \mu_n$,

$$f := \sum_{|\mu|=d} a_\mu X_0^{\mu_0} \cdots X_n^{\mu_n}, \qquad g := \sum_{|\mu|=d} a_\mu X_0^{\mu_0} \cdots X_n^{\mu_n}$$

and $\overline{b_\mu}$ stands for the complex conjugate of $b_\mu$.

We introduce a well–ordering in the set of multi–indices $S_d := \{(\mu_0, \ldots, \mu_n) : \mu_0 + \cdots + \mu_n = d\}$ as a bijection

$$\varphi : S_d \longrightarrow \{1 \leq i \leq N_d\}.$$

The reader may assume that $\varphi$ is given by the lexicographic order on $S_d$. Using this ordering we define the diagonal matrix $\Delta_d$ in the following terms :

$$\Delta_d := \left( \binom{d}{\mu_0 \cdots \mu_n}^{-1/2} \right)_{1 \leq \varphi(\mu_0, \ldots, \mu_n) \leq N_d},$$

where

$$\binom{d}{\mu_0 \cdots \mu_n} := \left( \frac{d!}{\mu_0! \cdots \mu_n!} \right).$$

As in [3], we also define an Hermitian inner product $< \cdot, \cdot >_d \; : H_d \times H_d \longrightarrow \mathbb{C}$ in the following terms.

$$< f, g >_d := < \Delta_d f, \Delta_d g > .$$

For every list $(d) := (d_1, \ldots, d_n) \in \mathbb{N}^n$ of degrees let $\mathcal{H}_{(d)} = H_{d_1} \times \cdots \times H_{d_n}$ be the space of sequences $F := (f_1, \ldots, f_n)$ of homogeneous polynomials such that $f_i \in H_{d_i}$, $1 \leq i \leq n$. We may also see $\mathcal{H}_{(d)}$ as the space of all polynomial mappings $F := (f_1, \ldots, f_n) \; : \mathbb{C}^{n+1} \longrightarrow \mathbb{C}^n$, such that $f_i \in H_{d_i}$, for all $i$.

We may extend the previous Hermitian inner products to the product space $\mathcal{H}_{(d)}$ in the following terms. Given $F := (f_1, \ldots, f_n) \in \mathcal{H}_{(d)}$ and $G := (g_1, \ldots, g_n) \in \mathcal{H}_{(d)}$, we define the *canonical Hermitian inner product* on $\mathcal{H}_{(d)}$ as $< \cdot, \cdot > \; : \mathcal{H}_{(d)} \times \mathcal{H}_{(d)} \longrightarrow \mathbb{R}$ in the following terms :

$$< F, G > := \sum_{i=1}^n < f_i, g_i > .$$

For every $F \in \mathcal{H}_{(d)}$, we denote by $\|F\|$ the norm defined by this canonical Hermitian inner product, i.e. $\|F\| := (< F, F >)^{\frac{1}{2}}$.

On the other hand, let $\Delta_{(d)}$ be the diagonal matrix given as the diagonal sum of $\Delta_{d_1}, \ldots, \Delta_{d_n}$. Namely,

$$\Delta_{(d)} := \Delta_{d_1} \oplus \cdots \oplus \Delta_{d_n}.$$

Then we also introduce the Hermitian inner product used in [3] in the following terms. Given $F, G \in \mathcal{H}_{(d)}$, we define $< F, G >_\Delta \in \mathbb{R}$ as

$$< F, G >_\Delta := < \Delta_{(d)} F, \Delta_{(d)} G > \in \mathbb{R}.$$

For every $F \in \mathcal{H}_{(d)}$, we denote by $\|F\|_\Delta$ the norm defined by this canonical Hermitian inner product, i.e. $\|F\|_\Delta := (< F, F >_\Delta)^{\frac{1}{2}}$.

Let $U(n+1)$ be the unitary group of isometries of $\mathbb{C}^{n+1}$. We may define the action of the group $U(n+1)$ on $\mathcal{H}_{(d)}$ in the following terms. Given $F \in \mathcal{H}_{(d)}$ and given $\sigma \in U(n+1)$, we define $\sigma(F) \in \mathcal{H}_{(d)}$ as the unique polynomial mapping in $\mathcal{H}_{(d)}$ that satisfies the following identity.

$$\sigma(F)(x) := F(\sigma^{-1}(x)), \; \forall x \in \mathbb{C}^{n+1}.$$

Then, the following statements holds

**Theorem 20** *[3] This Hermitian inner product $< \cdot, \cdot >_\Delta$ on $\mathcal{H}_{(d)}$ is unitarily invariant. In other words,*

$$< \sigma(F), \sigma(G) >_\Delta = < F, G >_\Delta,$$

*for every $\sigma \in U(n+1)$ and for every $F, G \in \mathcal{H}_{(d)}$.*

## 4.2 Two Riemannian structures on $\mathbb{P}(\mathcal{H}_{(d)})$.

Let $\mathbb{P}(\mathcal{H}_{(d)})$ be the complex projective space defined by the complex vector space $\mathcal{H}_{(d)}$. Hence, $N \in \mathbb{N}$ is the complex dimension of $\mathbb{P}(\mathcal{H}_{(d)})$. Recall that the real dimension of $\mathbb{P}(\mathcal{H}_{(d)})$ is $2N$ and that $N$ is given by the following identity :

$$N := \left( \sum_{i=1}^{n} N_{d_i} \right) - 1.$$

The well–order on the monomial exponents introduced in in Subsection 1.2 of the Introduction, yields an identification between $\mathcal{H}_{(d)}$ and $\mathbb{C}^{N+1}$. Accordingly, we may also identify $\mathbb{P}(\mathcal{H}_{(d)})$ and $\mathbb{P}_N(\mathbb{C})$.

We may introduce two Riemannian metrics in $\mathbb{P}_N(\mathbb{C})$ according to the Hermitian inner products on $\mathbb{C}^{N+1}$ discussed in Subsection 4.1 above. These two Riemannian metrics yield two Riemannian structures in $\mathbb{P}_N(\mathbb{C})$. We denote by $(\mathbb{P}_N(\mathbb{C}), can)$ the Riemannian structure defined by the canonical Hermitian inner product $< \cdot, \cdot > \; : \mathcal{H}_{(d)} \times \mathcal{H}_{(d)} \longrightarrow \mathbb{C}$. We denote by $(\mathbb{P}_N(\mathbb{C}), uni)$ the Riemannian structure defined by the (unitarily invariant) Hermitian inner product $< \cdot, \cdot >_\Delta \; : \mathcal{H}_{(d)} \times \mathcal{H}_{(d)} \longrightarrow \mathbb{C}$ introduced in Subsection 4.1 above.

Finally, let us denote by $\pi_\mathbb{C} \; : \mathbb{C}^{N+1} \setminus \{(0, \ldots, 0)\} \longrightarrow \mathbb{P}_N(\mathbb{C})$ the canonical projection given by

$$\pi_\mathbb{C}(z_0, \ldots, z_N) := (z_0 : z_1 : \ldots : z_N) \in \mathbb{P}(\mathcal{H}_{(d)}),$$

where $(z_0 : z_1 : \ldots : x_N)$ are the homogeneous coordinates of the projective point defined by $(z_0, \ldots, z_N) \in \mathbb{C}^{N+1}$.

Then the following Lemma holds :

**Lemma 21** *The following mapping is an isometry between $(\mathbb{P}(\mathcal{H}_{(d)}), can)$ and $(\mathbb{P}(\mathcal{H}_{(d)}), uni)$ :*

$$\widetilde{\Delta}^{-1} : (\mathbb{P}_N(\mathbb{C}), can) \longrightarrow (\mathbb{P}_N(\mathbb{C}), uni),$$

*where*
$$\widetilde{\Delta}^{-1}(z_0 : z_1 : \ldots : z_N) := \pi_{\mathbb{C}}(\Delta_{(d)}^{-1}(z_0, z_1, \ldots, z_N))),$$

*and $\Delta_{(d)}^{-1} \in GL(N+1, \mathbb{C})$ is the inverse of the regular matrix $\Delta_{(d)}$ introduced in Subsection 4.1 above.*

We denote by $\widetilde{\Delta}$ the inverse of the isometry $\widetilde{\Delta}^{-1}$ introduced above. Observe that for every $(z_0 : z_1 : \ldots : z_N) \in \mathbb{P}_N(\mathbb{C})$ the following holds :

$$\widetilde{\Delta}(z_0 : z_1 : \ldots : z_N) := \pi_{\mathbb{C}}(\Delta_{(d)}(z_0, z_1, \ldots, z_N))).$$

We may associate two volume forms to $\mathbb{P}_N(\mathbb{C})$ accordingly to the corresponding Riemannian structure. For every subset $B \subseteq \mathbb{P}_N(\mathbb{C})$ we denote as $Vol_{can}(B)$ the volume of $B$ with respect to the Riemannian structure $(\mathbb{P}_N(\mathbb{C}), can)$ and we denote as $Vol_{uni}(B)$ the corresponding volume of $B$ with respect to the Riemannian structure $(\mathbb{P}_N(\mathbb{C}), uni)$. From Lemma 21 above the following holds for every subset $B \subseteq \mathbb{P}_N(\mathbb{C})$ :

$$Vol_{can}(\widetilde{\Delta}(B)) = Vol_{uni}(B). \tag{18}$$

Now let $S^{2N+1} \subseteq \mathbb{R}^{2N+2}$ be the unit sphere. The unit sphere has also a canonical Riemannian structure defined by the canonical inner product in $\mathbb{R}^{2N+2}$. Identifying $\mathbb{R}^{2N+2} \cong \mathbb{C}^{N+1}$ we may also consider the canonical projection

$$p := \pi_{\mathbb{C}} \mid_{S^{2N+1}} : (S^{2N+1}, can) \longrightarrow (\mathbb{P}_N(\mathbb{C}), can).$$

The following Proposition is a well–known fact in Riemannian Geometry.

**Proposition 22** *The mapping $p$ is a Riemannian submersion. In fact, the following holds :*
$$S^{2N+1}/S^1 \cong \mathbb{P}_N(\mathbb{C}).$$

This Proposition yields the following identity for every subset $B \subseteq \mathbb{P}_N(\mathbb{C})$ :

$$Vol_S(p^{-1}(B)) = 2\pi Vol_{can}(B), \tag{19}$$

where $Vol_S(p^{-1}(B))$ denotes the volume of $p^{-1}(S)$ with respect to the canonical volume form (spherical volume) in $S^{2N+1}$.

Together with the Riemannian structures we have the corresponding *Fubiny–Study metrics* in $\mathbb{P}_N(\mathbb{C})$. We reproduce here this metric for the Riemannian structure $(\mathbb{P}_N(\mathbb{C}), uni)$. It is similarly defined in the other case. Recall that the Riemannian metric in $(\mathbb{P}_N(\mathbb{C}), uni)$ induces the distance function $d_R : \mathbb{P}_N(\mathbb{C}) \times \mathbb{P}_N(\mathbb{C}) \longrightarrow \mathbb{R}_+$ given by the following identity :

$$d_R(\pi_\mathbb{C}(X), \pi_\mathbb{C}(Y)) := \arccos\left(\frac{Re <X,Y>_\Delta}{\|X\|_\Delta \|Y\|_\Delta}\right),$$

where $Re <X,Y>_\Delta$ is the "real part" of the complex number $<X,Y>_\Delta \in \mathbb{C}$.

We define the Fubini–Study metric associated to the Riemannian structure $(\mathbb{P}_N(\mathbb{C}), uni)$ as the mapping $d_{FS} : \mathbb{P}_N(\mathbb{C}) \times \mathbb{P}_N(\mathbb{C}) \longrightarrow \mathbb{R}_+$ given by the following identity :

$$d_{FS}(x,y) := \sin \, d_R(x,y), \, \forall x, y \in \mathbb{P}_N(\mathbb{C}).$$

Observe that the Fubini–Study metric satisfies the following property.

**Proposition 23** *Let $B \subseteq \mathbb{P}_N(\mathbb{C}^N)$ be a subset and $\widetilde{B} := \pi_\mathbb{C}^{-1}(B) \cup \{0\}$. Let $X \in \mathbb{C}^{N+1}$ be a point such that*

$$\|X\|_\Delta = 1.$$

*Then, the following identity holds :*

$$inf\{\|X - Y\|_\Delta \, : \, Y \in \pi_\mathbb{C}^{-1}(B) \cup \{0\}\} = d_{FS}(\pi_\mathbb{C}(X), B)$$

*where $d_{FS}(\pi_\mathbb{C}(X), B) := inf\{d_{FS}(\pi_\mathbb{C}(X), y) \, : \, y \in B\}$.*

An analogous statement holds for the corresponding canonical Riemannian structure on $\mathbb{P}_N(\mathbb{C})$.

### 4.3 Solution and Discriminant Varieties

We may define the *solution variety*

$$V := \{(F, \zeta) \in \mathbb{P}(\mathcal{H}_{(d)}) \times \mathbb{P}(\mathbb{C}^{n+1}) : F(\zeta) = 0 \in \mathbb{C}^n\}.$$

**Proposition 24** *[3] With the previous notations and assumptions, $V$ is a smooth connected projective subvariety of $\mathbb{P}(\mathcal{H}_{(d)}) \times \mathbb{P}(\mathbb{C}^{n+1})$ of codimension $n$.*

Let us define the projection $\pi_1 : V \longrightarrow \mathbb{P}(\mathcal{H}_{(d)})$ and for every $(F, \zeta) \in V$, let

$$D\pi_1(F, \zeta) : T_{(F,\zeta)}V \longrightarrow T_F\mathbb{P}(\mathcal{H}_{(d)}),$$

the corresponding tangent mapping at $(F, \zeta)$. We say that $(F, \zeta) \in V$ is ill–conditioned if $D\pi_1(F, \zeta)$ is not an embedding. We define the set of ill–conditioned inputs $\Sigma' \subseteq V$ by the following identity

$$\Sigma' := \{(F, \zeta) \in V \, : \, \det D\pi_1(F, \zeta) = 0\}.$$

Then, we define the *discriminant variety* as the algebraic variety

$$\Sigma := \pi_1(\Sigma') := \pi_1(\{(F,\zeta) \in V : \det D\pi_1(F,\zeta) = 0\}) \subseteq \mathbb{P}(\mathcal{H}_{(d)}).$$

Let $\pi_2 : V \longrightarrow \mathbb{P}(\mathbb{C}^{n+1})$ be the projection onto the second factor. Namely $\pi_2(F,\zeta) := \zeta$ for all $(F,\zeta) \in V$.
Given $\zeta \in \mathbb{P}(\mathbb{C}^{n+1})$, we define $V_\zeta := \pi_2^{-1}(\zeta)$. Namely, $V_\zeta$ is the algebraic set of all systems of polynomial equations $F$ that vanishes at $\zeta$, i.e.

$$V_\zeta := \{(F,x) \in \mathbb{P}(\mathcal{H}_{(d)}) \times \mathbb{P}(\mathbb{C}^{n+1}) : (F,x) \in V \wedge x = \zeta\}.$$

Given $(F,\zeta) \in V$, we define the *fiber distance* to the discriminant variety $\Sigma$ in the following terms :

$$\rho(F,\zeta) := \inf\{d_{FS}(F,G) \; : \; (G,\zeta) \in \Sigma' \bigcap V\}.$$

**Lemma 25** *[3] The function $\rho \; : V \longrightarrow \mathbb{R}$ is invariant under the action of $U(n+1)$ on $V$. Namely, given $(F,\zeta) \in V$ and given $\sigma \in U(n+1)$, the following equality holds :*

$$\rho(\sigma(F), \sigma(\zeta)) = \rho(F,\zeta).$$

We also define the fiber distance of system $F \in V \subseteq \mathbb{P}(\mathcal{H}_{(d)})$ to the discriminant variety $\Sigma$ as

$$\rho(F) := \min\{\rho(F,\zeta) \; : \; F(\zeta) = 0\}.$$

Obviously, $\rho(F) = 0$ when $F \in \Sigma$.
On the other hand, Let $F := (f_1, \ldots, f_n) \in \mathcal{H}_{(d)}$ be a system of homogeneous polynomials. Let $F \; : \mathbb{C}^{n+1} \longrightarrow \mathbb{C}^n$ be the regular mapping between affine spaces defined by $F$. For every $z \in \mathbb{C}^{n+1}$, let $DF(z) \; : \mathbb{C}^{n+1} \longrightarrow \mathbb{C}^n$ be the tangent mapping defined by the jacobian of $F$ at $z$. For every $\zeta \in \mathbb{C}^{n+1}$, let $T_\zeta := \{w \in \mathbb{C}^{n+1} \; : \; <w,\zeta> = 0\}$ the orthogonal complement of $\zeta \in \mathbb{C}^{n+1}$ with respect to the canonical Hermitian inner product.
Let us denote by $Diag(\|\zeta\|^{d_i-1} d_i^{1/2})$ the square diagonal matrix given by the following identity :

$$Diag(\|\zeta\|^{d_i-1} d_i^{1/2}) := \begin{pmatrix} \|\zeta\|^{d_1-1} d_1^{1/2} & \cdots & 0 \\ \vdots & \ddots & \vdots \\ 0 & \cdots & \|\zeta\|^{d_n-1} d_n^{1/2} \end{pmatrix}.$$

Let $F \in \mathcal{H}_{(d)}$ be a system of polynomial equations and let $\zeta \in \mathbb{P}_n(\mathbb{C})$ be a non–singular zero (i.e. $F(\zeta) = 0$ and the linear mapping $DF(\zeta) : T_\zeta \longrightarrow \mathbb{C}^n$ is a non–singular linear mapping). The normalized condition number of $F$ at $\zeta$ is define by the following identity :

$$\mu_{norm}(F,\zeta) := \|F\|_\Delta \|DF(\zeta) \mid_{T_\zeta}^{-1} Diag(\|\zeta\|^{d_i-1} d_i^{1/2})\|,$$

where $\|F\|_\Delta$ denotes the norm of $F$ with respect to the unitarily invariant Hermitian inner product in $\mathcal{H}_{(d)}$ introduced in Subsection 4.1 above. We define the normalized condition number $\mu_{norm}(F)$ as the maximum :

$$\mu_{norm}(F) := \max\{\mu_{norm}(F, \zeta) : \zeta \in \mathbb{P}_n(\mathbb{C}), F(\zeta) = 0\}.$$

In these pages, we are not going to discuss the relevance of this notion both for the behaviour of the projective Newton's method and for the complexity of homotopy continuation methods to solve systems of multivariate polynomial equations. The reader interested may follow the comprehensive discussion in Chapter 14 of [3] and the references therein. The following central statement explains one of the main features of this normalized condition number.

**Theorem 26 (Condition Number Theorem, [3])** *With the previous notations and assumptions, for every system $F \in \mathbb{P}(\mathcal{H}_{(d)})$, the following holds :*

$$\mu_{norm}(F) = \frac{1}{\rho(F)}.$$

In fact, the Condition Number Theorem above gives a more precise statement that we are going to use in the sequel. This is explicitly given in the following statement.

**Proposition 27** *[3] Let $F := (f_1, \ldots, f_n) \in \mathcal{H}_{(d)}$ be a system of homogeneous polynomials and $\zeta \in \mathbb{P}_n(\mathbb{C})$. Assume $(F, \zeta) \in V$ and $\|F\|_\Delta = 1$. Let $\sigma \in U(n+1)$ be a unitary transformation such that $\sigma(\zeta) := (1, 0, \ldots, 0) = e_0 \in \mathbb{C}^{n+1}$. Then, the following identity holds :*

$$\rho(F, \zeta) := d_F(D\sigma(F)(e_0) Diag(d_i^{-1/2}), S),$$

*where*

   i) $S \subseteq \mathcal{M}_n(\mathbb{C})$ *is the algebraic variety of singular matrices,*

   ii) $D\sigma(F)(e_0) \in \mathcal{M}_n(\mathbb{C})$ *is the jacobian matrix :*

$$D\sigma(F)(e_0) := \left(\frac{\partial \sigma(f_i)}{\partial X_j}(e_0)\right)_{1 \leq i,j \leq n},$$

   iii) $Diag(d_i^{-1/2})$ *is the diagonal matrix given by the following identity :*

$$Diag(d_i^{-1/2}) := \begin{pmatrix} d_1^{-1/2} & \cdots & 0 \\ \vdots & & \vdots \\ 0 & \cdots & d_n^{-1/2} \end{pmatrix}.$$

iv) $d_F$ is the Frobenius distance introduced in Subsection 2.2 above, namely given $X := (x_{i,j})$ and $Y := (y_{i,j}) \in \mathcal{M}_n(\mathbb{C})$, we define

$$d_F(X,Y) := \|X - Y\|_F := \left( \sum_{1 \leq i,j \leq} |x_{i,j} - y_{i,j}|^2 \right)^{1/2}.$$

For every $\varepsilon > 0$ we define the following tubular neighborhood of the discriminant variety $\Sigma$.

$$\Sigma(\varepsilon) := \{F \in V \subseteq \mathbb{P}(\mathcal{H}_{(d)}) : \rho(F) \leq \varepsilon\}.$$

Now, the following statement holds.

**Theorem 28** *[3] With the previous notations, the following holds :*

$$\frac{Vol_{uni}(\Sigma(\varepsilon))}{Vol_{uni}(\mathbb{P}(\mathcal{H}_{(d)}))} \leq \varepsilon^4 \mathfrak{C}[(d)],$$

*where $\mathfrak{C}[(d)] := n^3(n+1)N(N-1)\mathcal{D}_{(d)}$ is the constant introduced in Theorem 3 above.*

For every positive real number $H > 0$, let $B_\Delta(0, H) \subseteq \mathbb{C}^{N+1}$ be the closed ball of radius $H$ centered at the origin $0 \in \mathbb{C}^{N+1}$ with respect to the unitarily invariant norm $\|\cdot\|_\Delta$. Namely,

$$B_\Delta(0, H) := \{Z \in \mathbb{C}^{N+1} : \|Z\|_\Delta \leq H\}.$$

On the other hand, let $B(0, H) \subseteq \mathbb{C}^{N+1}$ the closed ball of radius $H$ centered at the origin $0 \in \mathbb{C}^{N+1}$ with respect to the canonical norm $\|\cdot\|$. The following equality holds :

$$B_\Delta(0, H) = \Delta_{(d)}^{-1} B(0, H), \tag{20}$$

where $\Delta_{(d)} : \mathbb{C}^{N+1} \longrightarrow \mathbb{C}^{N+1}$ is the complex vector space automorphism given by the regular matrix $\Delta_{(d)}$.

With the same notations as above, for every $\varepsilon > 0$ and for every positive real number $H > 0$, let $\overline{R}_\Delta(\varepsilon, H)$ be the compact set given by the following identity

$$\overline{R}_\Delta(\varepsilon, H) := \left( \pi_\mathbb{C}^{-1}(\Sigma(\varepsilon)) \bigcup \{0\} \right) \bigcap B_\Delta(0, H). \tag{21}$$

The following statement is a consequence of Theorem 28 above :

**Corollary 29** *With the previous notations and assumptions, the Lebesgue measure of $\overline{R}_\Delta(\varepsilon, H)$ satisfies the following inequality :*

$$Vol(\overline{R}_\Delta(\varepsilon, H)) \leq \varepsilon^4 H^{2N+2} \mathfrak{D}[(d)],$$

*where*

$$\mathfrak{D}[(d)] := \frac{2\pi Vol_{uni}(\mathbb{P}(\mathcal{H}_{(d)}))\mathfrak{C}[(d)]}{(2N+2)\det(\Delta_{(d)})} = \frac{K_{2N+2}}{\det(\Delta_{(d)})} \mathfrak{C}[(d)].$$

*Proof.–* For every $H$, let us define the compact subset $V(\varepsilon, H)$ given by the following identity :

$$V(\varepsilon, H) := \left(\pi_{\mathbb{C}}^{-1}(\widetilde{\Delta}(\Sigma(\varepsilon))) \bigcup \{0\}\right) \bigcap B(0, H),$$

where

$$\widetilde{\Delta}(\Sigma(\varepsilon))) := \{\widetilde{\Delta}(x) \ : \ x \in \Sigma(\varepsilon)\}.$$

By integration in spherical coordinates we obtain

$$Vol(V(\varepsilon, H)) = \frac{H^{2N+2}}{2N+2} Vol_S(p^{-1}(\widetilde{\Delta}(\Sigma(\varepsilon)))).$$

Now by Identity (19) above we conclude

$$Vol(V(\varepsilon, H)) = \frac{2\pi H^{2N+2}}{2N+2} Vol_{can}(\widetilde{\Delta}(\Sigma(\varepsilon)))).$$

Finally, from Identity (18) above, we conclude :

$$Vol(V(\varepsilon, H)) = \frac{2\pi H^{2N+2}}{2N+2} Vol_{uni}(\Sigma(\varepsilon)).$$

On the other hand, the following identity holds :

$$\pi_{\mathbb{C}}^{-1}(\widetilde{\Delta}(\Sigma(\varepsilon))) = \Delta_{(d)}(\pi_{\mathbb{C}}^{-1}(\Sigma(\varepsilon))).$$

As $B_\Delta(0, H) = \Delta_{(d)}^{-1} B(0, H)$, we conclude :

$$V(\varepsilon, H) = \Delta_{(d)} R_\Delta(0, H),$$

Thus, we conclude

$$Vol(\overline{R}_\Delta(\varepsilon, H)) \leq \frac{Vol(V(\varepsilon, H))}{\det(\Delta_{(d)})} = \frac{2\pi H^{2N+2}}{(2N+2)\det(\Delta_{(d)})} Vol_{uni}(\Sigma(\varepsilon)),$$

and the Proposition follows from Theorem 28 above. ∎

## 4.4 Gauss integers, unitarily invariant height and $\mathbb{C}-$visible points

The theory of visible points (as used in previous pages) is not directly applicable to the study of points in a complex projective space. In this Subsection we shall see how to modify that theory of visible points to count points in a projective space $\mathbb{P}_m(\mathbb{Q}[i])$, where $\mathbb{Q}[i]$ is the field of Gauss rationals.
Let $\mathbb{Z}[i] = \{a + bi \ : \ a, b \in \mathbb{Z}\}$ be the ring of Gauss integers and let $\mathbb{Z}[i]^*$ be the group of units. Namely $\mathbb{Z}[i]^* := \{1, -1, i, -i\}$, where $i^2 = -1$. As $\mathbb{Z}[i]$ is a factorial domain, we may also introduce the corresponding theory of visible points in a $\mathbb{Z}[i]-$lattice (compare with Subsection 2.4 above).

A $\mathbb{Z}[i]-$lattice in $\mathbb{C}^{m+1}$ is the free $\mathbb{Z}[i]-$module generated by a basis of $\mathbb{C}^{m+1}$ as complex vector space. For instance, a $\mathbb{Z}[i]-$lattice is the set $\mathbb{Z}[i]^{m+1} \subseteq \mathbb{C}^{m+1}$ of all points in $\mathbb{C}^{m+1}$ whose coordinates are Gauss integers.

In these pages we consider the following two examples of $\mathbb{Z}[i]-$ lattices. First of all, the $\mathbb{Z}[i]-$lattice $\mathcal{H}_{(d)}(\mathbb{Z}[i])$ of all sequences of homogeneous polynomials $F := (f_1, \ldots, f_n) \in \mathcal{H}_{(d)}$ whose coefficients are Gauss integers (namely, $f_i \in \mathbb{Z}[i][X_0, \ldots, X_n]$). For every complex point $\rho \in \mathbb{C}$, we also consider the $\mathbb{Z}[i]-$lattice $\mathcal{H}_{(d)}(\mathbb{Z}[i])(\rho)$ given by the following identity :

$$\mathcal{H}_{(d)}(\mathbb{Z}[i])(\rho) := \{\rho F \ : \ F \in \mathcal{H}_{(d)}(\mathbb{Z}[i])\}.$$

**Definition 30** *A non–zero point $\rho F := (\rho z_0, \ldots, \rho z_N) \in \mathcal{H}_{(d)}(\mathbb{Z}[i])(\rho)$ is said to be $\mathbb{C}-$visible from the origin if and only if the following holds :*

$$gcd_{\mathbb{Z}[i]}(F) := gcd_{\mathbb{Z}[i]}\{z_0, \ldots, z_N\} \in \mathbb{Z}[i]^*,$$

*where $gcd_{\mathbb{Z}[i]}\{z_0, \ldots, z_N\}$ is the greatest common divisor of $z_0, \ldots, z_N$ in $\mathbb{Z}[i]$.*

Observe that for $\rho \in \mathbb{R} \setminus \{0\}$, the class of $\mathbb{C}$–visible points of $\mathcal{H}_{(d)}(\mathbb{Z}[i])(\rho)$ is not equal to the class of visible points of $\mathbb{Z}^{2N+2}(\rho)$ (in the sense of Subsection 2.4 above).

Let $\mathbb{P}(\mathcal{H}_{(d)})$ be the complex projective space of dimension $N$ defined above and let $\mathbb{P}(\mathcal{H}_{(d)}(\mathbb{Q}[i]))$ be the $N-$dimensional projective space defined by the $\mathcal{H}_{(d)}(\mathbb{Q}[i])$. Let $\pi_{\mathbb{C}} : \mathcal{H}_{(d)} \setminus \{0\} \longrightarrow \mathbb{P}(\mathcal{H}_{(d)})$ be the canonical projection. Let us also denote by $\pi_{\mathbb{C}}$ the canonical projection

$$\pi_{\mathbb{C}} : \mathcal{H}_{(d)}(\mathbb{Q}[i]) \setminus \{0\} \longrightarrow \mathbb{P}(\mathcal{H}_{(d)}(\mathbb{Q}[i])).$$

In order to introduce the notion of *unitarily invariant height of a projective point in $\mathbb{P}(\mathcal{H}_{(d)}(\mathbb{Q}[i]))$* (see also Subsection 2.4 for a comparison) we resume here in a very concise form the language and notation used for absolute values over number fields. For an introduction refer to e.g. [17, Chapter 1], whereas a more complete exposition of the theory of absolute values can be found in Artin's *Algebraic Numbers and Algebraic Functions* [1] or [18].

Let $|\cdot|_\nu : \mathbb{Q}[i] \longrightarrow \mathbb{R}_+$ be an absolute value defined on the number field $\mathbb{Q}[i]$. By $\mathbb{Q}[i]_\nu$ we denote the completion of $\mathbb{Q}[i]$ with respect to this absolute value $|\cdot|_\nu$ and by $\overline{\mathbb{Q}[i]}_\nu$ we denote the algebraic closure of $\mathbb{Q}[i]_\nu$. Finally, we denote by $n_\nu$ the degree of $\mathbb{Q}[i]_\nu$ over the completion of $\mathbb{Q}$ with respect to the absolute value $|\cdot|_\nu : \mathbb{Q} \longrightarrow \mathbb{R}_+$.

Let $M_{\mathbb{Q}[i]}$ be a proper set of absolute values of $\mathbb{Q}[i]$ in the sense of [17]. We assume that $M_{\mathbb{Q}[i]}$ has been chosen such that it satisfies Weil's *product formula* with multiplicities $n_\nu$, i.e. for all $z \in \mathbb{Q}[i] \setminus \{0\}$ the following holds

$$\prod_{\nu \in M_{\mathbb{Q}[i]}} |z|_\nu^{n_\nu/2} = 1 \qquad (22)$$

Let $\|\cdot\| : \mathbb{C}^{m+1} \longrightarrow \mathbb{R}_+$ the canonical Hermitian norm in $\mathbb{C}^{m+1}$. Let $S \subseteq M_{\mathbb{Q}[i]}$ be the class of all sub–indices $\nu \in M_{\mathbb{Q}[i]}$ such that the absolute value $|\cdot|_\nu : \mathbb{Q}[i] \longrightarrow \mathbb{R}_+$ is a non–archimedean absolute value.

**Definition 31** *For every point $\pi_\mathbb{C}(F) := (z_0 : \ldots : z_N) \in \mathbb{P}(\mathcal{H}_{(d)}(\mathbb{Q}[i]))$, where $z_0, \ldots, z_N \in \mathbb{Q}[i]$, we define unitarily invariant height of $\pi_\mathbb{C}(F)$ as :*

$$H_\Delta(\pi_\mathbb{C}(F)) := \|F\|_\Delta \left( \prod_{\nu \in S} \max\{|z_0|_\nu, \ldots, |z_N|_\nu\}^{n_\nu} \right)^{1/2}.$$

Since Weil's product formula (Equation (22) above) this quantity is well-defined and independent of the affine point $F = (z_0, \ldots, z_N) \in \mathcal{H}_{(d)}(\mathbb{Q}[i]) \setminus \{0\}$ we have chosen.

**Remark 32** *We may also define the Northcott–Schmidt height of a projective point $\pi_\mathbb{C}(F) \in \mathbb{P}(\mathcal{H}_{(d)}(\mathbb{Q}[i]))$ as in Subsection 2.4 above. Namely, we may introduce*

$$H(\pi_\mathbb{C}(F)) := \|F\| \left( \prod_{\nu \in S} \max\{|z_0|_\nu, \ldots, |z_N|_\nu\}^{n_\nu} \right)^{1/2}.$$

*Let us observe that this notion of height is essentially equivalent to Weil's absolute height of the projective point $\pi_\mathbb{C}(z)$ as used in [13], [16] or [4] and the references therein.*
*From the Definition of $\Delta_{(d)}$ one easily concludes :*

$$\frac{H(\pi_\mathbb{C}(F))}{D_{(d)}!} \leq H_\Delta(\pi_\mathbb{C}(F)) \leq H(\pi_\mathbb{C}(F)). \tag{23}$$

*However unitarily invariant height is better suited for our purposes that the original of Northcott and Schmidt.*

Unitarily invariant height and $\mathbb{C}$–visible points are closely related.

**Lemma 33** *With the same notations and assumptions as above, the following properties hold :*

i) *Let $\rho F \in \mathcal{H}_{(d)}(\mathbb{Z}[i])(\rho)$ be a point in that $\mathbb{Z}[i]$–lattice. Then, $\rho F$ is a $\mathbb{C}$–visible point in $\mathcal{H}_{(d)}(\mathbb{Z}[i])(\rho)$ if and only if*

$$\|F\|_\Delta = min\{\|G\|_\Delta : \rho G \in \mathcal{H}_{(d)}(\mathbb{Z}[i])(\rho) \setminus \{0\}, \pi_\mathbb{C}(\rho F) = \pi_\mathbb{C}(\rho G)\}.$$

ii) *Let $\rho F \in \mathcal{H}_{(d)}(\mathbb{Z}[i])(\rho)$ be a $\mathbb{C}$–visible point. Then,*

$$\|F\|_\Delta = H_\Delta(\pi_\mathbb{C}(\rho F)).$$

*iii)* Let $V \subseteq \mathbb{P}(\mathcal{H}_{(d)})$ be a subset and let $\widetilde{V} := \pi_{\mathbb{C}}^{-1}(V) \cup \{0\}$ be the complex cone over $V$. Let $B_\Delta(0, H)$ be the closed ball in $\mathbb{C}^{m+1}$ of radius $H$ centered at the origin. Then, the following holds :

- For every $\pi_{\mathbb{C}}(F) \in V \cap \mathbb{P}_m(\mathbb{Q}[i])$ there is a $\mathbb{C}$–visible point $G \in \widetilde{V} \cap \mathcal{H}_{(d)}(\mathbb{Z}[i])$ such that $\pi_{\mathbb{C}}(G) = \pi_{\mathbb{C}}(F)$ and $\|G\|_\Delta = H_\Delta(\pi_{\mathbb{C}}(F))$.
- The number of $\mathbb{C}$–visible points of $\mathcal{H}_{(d)}(\mathbb{Z}[i])$ that lie in $\widetilde{V} \cap B_\Delta(0, H)$ is four times the number of projective points in $V \cap \mathbb{P}(\mathcal{H}_{(d)}(\mathbb{Q}[i]))$ of unitarily invariant height at most $H$.

In fact, given projective point $\pi_{\mathbb{C}}(F) \in \mathbb{P}(\mathcal{H}_{(d)}(\mathbb{Q}[i]))$ and given a $\mathbb{C}$–visible point $G \in \mathcal{H}_{(d)}(\mathbb{Z}[i])$ such that $\pi_{\mathbb{C}}(G) = \pi_{\mathbb{C}}(F)$, the four $\mathbb{C}$–visible points in $\mathcal{H}_{(d)}(\mathbb{Z}[i])$ associated to $\pi_{\mathbb{C}}(F)$ are $\{G, -G, iG, -iG\} \subseteq \mathcal{H}_{(d)}(\mathbb{Z}[i])$.
Given a projective point $\pi_{\mathbb{C}}(F) \in \mathbb{P}(\mathcal{H}_{(d)}(\mathbb{Q}[i]))$, we may represent this point by means of some $\mathbb{C}$–visible point $G \in \mathcal{H}_{(d)}(\mathbb{Z}[i])$ such that $\pi_{\mathbb{C}}(G) = \pi_{\mathbb{C}}(F)$. Thus, we may define the bit length of the projective point $\pi_{\mathbb{C}}(F)$ as the bit length of this representation $G$. Actually, for points in $\mathcal{H}_{(d)}(\mathbb{Z}[i])$, the number of tape cells required to represent $G$ in a Turing machine essentially agrees with the logarithm of its norm, i.e. $\log \|G\|$. As in Subsection 2.4 above, we may define the bit length of $\pi_{\mathbb{C}}(F)$ as the logarithm of its Northcott–Schmidt height. Namely

$$bl(\pi_{\mathbb{C}}(F)) := \log \|G\| = \log H(\pi_{\mathbb{C}}(F)).$$

As observed in Remark 32 above, Northcott–Schmidt height and unitarily invariant height are essentially equivalent. In Fact, the following holds for every $\pi_{\mathbb{C}}(F) \in \mathbb{P}(\mathcal{H}_{(d)}(\mathbb{Q}[i]))$ :

$$bl(\pi_{\mathbb{C}}(F)) - D_{(d)} \log D_{(d)} \leq \log H_\Delta(\pi_{\mathbb{C}}(F)) \leq bl(F),$$

where $D_{(d)} := \max\{d_1, \ldots, d_n\}$. Thus we introduce the following notion :

**Definition 34** *For every $\pi_{\mathbb{C}}(F) \in \mathbb{P}(\mathcal{H}_{(d)}(\mathbb{Q}[i]))$ we define its unitarily invariant bit length as the logarithm of its unitarily invariant height. Namely,*

$$bl_\Delta(\pi_{\mathbb{C}}(F)) := \log H_\Delta(\pi_{\mathbb{C}}(F)).$$

## 5 The Non–linear Case

### 5.1 Davenport's discrepancy in the Non–linear Case

In this Subsection we apply Davenport's strategy of Subsection 2.3 above to the Non–linear case. Given two positive real numbers $H, \varepsilon \in \mathbb{R}$, let us denote by $\overline{R}_\Delta(\varepsilon, H)$ the compact subset defined in Equation 21 above. Namely,

$$\overline{R}_\Delta(\varepsilon, H) := \left(\pi_{\mathbb{C}}^{-1}(\Sigma(\varepsilon)) \bigcup \{0\}\right) \bigcap B_\Delta(0, H).$$

**Proposition 35** *With the same notations and assumptions as in Section 4 above, let $H, \varepsilon \in \mathbb{R}$ be two positive real numbers. Assume $H \geq 1$. Let $(d) := (d_1, \ldots, d_n)$ be a list of degrees and $D_{(d)} := \max\{d_1, \ldots, d_n\}$. Let $\overline{N}(\varepsilon, H)$ be the number of points of the $\mathbb{Z}[i]-$lattice $\mathcal{H}_{(d)}(\mathbb{Z}[i])$ that belong to $\overline{R}_\Delta(\varepsilon, H)$. Namely,*

$$\overline{N}(\varepsilon, H) := \sharp \left( \overline{R}_\Delta(\varepsilon, H) \cap \mathcal{H}_{(d)}(\mathbb{Z}[i]) \right).$$

*Then, the following inequality holds :*

$$|\overline{N}(\varepsilon, H) - Vol(\overline{R}_\Delta(\varepsilon, H))| \leq \frac{\overline{T}_N \mathfrak{S}^{(2N+2)}}{\det(\Delta_{(d)})^2} H^{2N+1},$$

*where $\mathfrak{S}^{(2N+2)}$ is the constant defined in Identity 1 of the Introduction and*

$$\overline{T}_N := \max\{8, 4(D_{(d)} + 1)\}^{4\left[N^2 + 3(n+2)^2(N+1)\right]}.$$

*Proof.–* First of all, observe that $\overline{R}_\Delta(\varepsilon, H)$ is the set of all systems of homogeneous polynomials $F := (f_1, \ldots, f_n) \in \mathcal{H}_{(d)}$ that satisfies the following properties :

- $\|F\|_\Delta \leq H$,
- $\exists \zeta \in \mathbb{C}^{n+1}$, such that $f_1(\zeta) = \cdots = f_n(\zeta) = 0$,
- $\rho(F, \zeta) \leq \varepsilon$.

From Proposition 27 above we may rewrite these three properties in the following terms. A polynomial $F \in \overline{R}(\varepsilon, H)$ if and only if the following properties hold :

- $\|F\|_\Delta^2 \leq H^2$,
- $\exists \zeta \in \mathbb{C}^{n+1}$, such that $f_1(\zeta) = \cdots = f_n(\zeta) = 0$,
- $\exists \sigma \in U(n+1)$, such that the following holds :
  - $\sigma(\zeta) = e_0 = (1, 0, \ldots, 0)$,
  - $\exists N \in \mathcal{M}_n(\mathbb{C})$ such that $\det(N) = 0$ and the following inequality holds
  
  $$\left\| D\sigma(F)(e_0) Diag(d_i^{-1/2}) - \|F\|_\Delta N \right\|_F^2 \leq \|F\|_\Delta^2 \varepsilon^2.$$

From now on, let us identify $\mathbb{C} \cong \mathbb{R}^2$, $\mathbb{Z}[i] \cong \mathbb{Z}^2$ and, accordingly, let us identify

$$\mathcal{M}_n(\mathbb{C}) \cong \mathbb{R}^{2n^2}, \ \mathcal{H}_{(d)} \cong \mathbb{R}^{2N+2}, \mathbb{C}^{n+1} \cong \mathbb{R}^{2n+2}.$$

We also denote by $U(n+1)$ the corresponding real algebraic variety in $\mathbb{R}^{2n^2}$. Observe that $U(n+1) \subseteq \mathbb{R}^{2n^2}$ is given by a finite number of quadratic equations (i.e. homogeneous equations of degree 2) that may be rewritten as $\sigma \overline{\sigma}^t = Id_{n+1}$.

Observe that for every $\sigma \in U(n+1)$, the inverse of $\sigma$ is the transpose conjugate of $\sigma$ (i.e. $\sigma^{-1} := \overline{\sigma}^t$). In particular, given $\sigma := (u_{i,j}) \in U(n+1)$ and given $F := (f_1, \ldots, f_n) \in \mathcal{H}_{(d)}$, the coefficients of $\sigma(f_i)$ are polynomials in the coefficients of $f_i$ and the entries $u_{i,j}$ of $\sigma$ of degree $d_i + 1$. This is going to be used in the sequel.

For every $m$, $1 \leq m \leq 2N+2$, let $\mathcal{I}_m$ be the class of all subsets $I \subseteq \{1, \ldots, 2N+2\}$ of cardinality $m$. For $I \in \mathcal{I}_m$, let $J \subseteq I$ be a subset such that $\sharp(I \setminus J) = 1$. Let $\alpha := (\alpha_i \ : \ j \in J) \in \mathbb{R}^{m-1}$ be an affine real point. Finally, let us consider the real algebraic variety :

$$V(\varepsilon, H, \alpha, I, J) \subseteq \mathbb{R}^{2N+2+2(n+1)+4(n+1)^2+4},$$

given by the following properties :
A point $(F, \zeta, \sigma, N, t_1, t_2, t_3, t_4) \in \mathbb{R}^{2N+2} \times \mathbb{R}^{2(n+1)} \times \mathbb{R}^{2(n+1)^2} \times \mathbb{R}^{2(n+1)^2} \times \mathbb{R}^4$
belongs to $V(\varepsilon, H, \alpha, I, J)$ if and only if the following properties hold :

i) $\|F\|_\Delta^2 - H = t_1^2$,

ii) $F_j = \alpha_j, \ \forall j \in J$, where $F := (F_1, \ldots, F_{2N+2}) \in \mathbb{R}^{2N+2}$.

iii) $F(\zeta) = (0, \ldots, 0)$, where $F$ is decomposed in the list of real and imaginary parts of the coefficients of the polynomials $f_i \in \mathcal{H}_{d_i}$, $1 \leq i \leq n$ such that $F := (f_1, \ldots, f_n)$, $\zeta := (\zeta_1, \ldots, \zeta_{2(n+1)}) \in \mathbb{R}^{2(n+1)}$ is identified with the corresponding complex point $\zeta \in \mathbb{C}^{n+1}$ and $F(\zeta)$ is understood as the real and imaginary parts of the complex numbers $F(\zeta) \in \mathbb{C}^n$.

iv) $\sigma(\zeta) = (1, 0, \ldots, 0) \in \mathbb{R}^{2(n+1)}$. Taking real and imaginary parts again these are $2(n+1)$ quadratic equations.

v) $\sigma \overline{\sigma}^t = Id_{n+1}$, taking real and imaginary parts these are also $2(n+1)$ quadratic equations.

vi) $\|F\|_\Delta^2 = t_3^2$ and $t_3 - t_4^2 = 0$. In other words, we rewrite as quadratic polynomial equations the fact $t_3$ is the norm of $F$ with respect to $<\cdot,\cdot>_\Delta$, i.e. $t_3 = \|F\|_\Delta$.

vii) and, finally

$$\left\|D\sigma(F)(e_0)Diag(d_i^{-1/2}) - t_3 N\right\|_F^2 - t_3^2 \varepsilon^2 = t_2^2. \qquad (24)$$

Observe that Equation (24) is a polynomial equation of degree at most $\max\{4, 2(D_{(d)} + 1)\}$, where $D_{(d)} := \max\{d_1, \ldots, d_n\}$.

Now, we are in conditions to apply the upper bounds for the sum of the Betti numbers of a real algebraic set, simultaneous and independently obtained by J. Milnor (cf. [19]), R. Thom (cf. [35]), O. Oleinik (cf. [23]) and A. Petrovsky (cf. [24]). These upper bounds imply the following inequality :

$$\beta_0(V(\varepsilon, H, \alpha, I, J)) \leq \big(\max\{8, 4(D_{(d)} + 1)\}\big)^{2N+2+2(n+1)+4(n+1)^2+4},$$

where $\beta_0(V(\varepsilon, H, \alpha, I, J))$ is the number of connected components of the real algebraic set $V(\varepsilon, H, \alpha, I, J)$.

Let $\Pi : \mathbb{R}^{2N+2+2(n+1)+4(n+1)^2+4} \longrightarrow \mathbb{R}^{2N+2}$ be the projection given by :

$$\Pi(F, \zeta, \sigma, N, t_1, t_2, t_3, t_4) := F \in \mathbb{R}^{2N+2}.$$

Next, let $\Phi_I : \mathbb{R}^{2N+2} \longrightarrow \mathbb{R}^m$ be the projection given by

$$\Phi_I(F) := (F_i \ : \ i \in I) \in \mathbb{R}^m.$$

Let $r_J(\alpha) \subseteq \mathbb{R}^m$ the real line given by the conditions :

$$r_J(\alpha) := \{(x_i \ : \ i \in I) \in \mathbb{R}^m \ : \ x_j = \alpha_j, \forall j \in J\}.$$

Then, the following holds :

$$\Phi_I(\Pi(V(\varepsilon, H, \alpha, I, J))) = \Phi_I(\overline{R}_\Delta(\varepsilon, H)) \bigcap r_J(\alpha).$$

As $\Phi_I$ and $\Pi$ are continuous mappings, we conclude

$$\beta_0(\Phi_I(\overline{R}_\Delta(\varepsilon, H) \bigcap r_J(\alpha))) \leq \big(\max\{8, 4(D_{(d)} + 1)\}\big)^{2N+6(n+1)^2+6}.$$

With the same notations as in Subsection 2.3 we conclude

$$h(\overline{R}_\Delta(\varepsilon, H), m) \leq \big(\max\{8, 4(D_{(d)} + 1)\}\big)^{2N+6(n+1)^2+6}. \tag{25}$$

On the other hand, for every $I \in \mathcal{I}_m$, the following holds

$$\Phi_I(\overline{R}_\Delta(\varepsilon, H)) \subseteq B_\Delta^{(m)}(0, H),$$

where $B_\Delta^{(m)}(0, H)$ is the closed ball in $\mathbb{R}^m$ of radius $H$, centered at the origin with the norm induced by $\|\cdot\|_\Delta$ on $\mathbb{R}^m \subseteq \mathbb{R}^{2N+2}$. Hence, we conclude from the Definition of $\Delta_{(d)}$ that the following holds :

$$Vol(\Phi_I(\overline{R}_\Delta(\varepsilon, H))) \leq H^m \frac{K_m}{\det(\Delta_{(d)})^2}.$$

With the notations of Subsection 2.3 the volume estimates $V(\overline{R}_\Delta(\varepsilon, H), m)$ satisfy the following estimate :

$$V(\overline{R}_\Delta(\varepsilon, H), m) := \sum_{I \in \mathcal{I}_m} Vol(\Phi_I(\overline{R}_\Delta(\varepsilon, H))) \leq \frac{H^m \binom{2N+2}{m} K_m}{\det(\Delta_{(d)})^2}. \quad (26)$$

The following inequality is a consequence of Theorem 12 :

$$|\overline{N}(\varepsilon, H) - Vol(\overline{R}_\Delta(\varepsilon, H)| \leq \sum_{\ell=0}^{2N+1} h(\overline{R}_\Delta(\varepsilon, H))^{2N+2-\ell} V(\overline{R}_\Delta(\varepsilon, H), \ell).$$

As $H \geq 1$, from Inequalities (25) and (26) we conclude :

$$|\overline{N}(\varepsilon, H) - Vol(\overline{R}_\Delta(\varepsilon, H)| \leq \frac{\overline{T}_N \mathfrak{S}^{(2N+2)}}{\det(\Delta_{(d)})^2} H^{2N+1},$$

where $\overline{T}_N$ and $\mathfrak{S}^{(2N+2)}$ are the constants stated above. ∎

The following Proposition follows by the same argument.

**Proposition 36** *With the same notations and assumptions as in Proposition 35 above, let $H \in \mathbb{R}$ be a positive real number, $H \geq 1$, and let $\overline{N}(1, H)$ be the number of systems of homogeneous polynomials in $\mathcal{H}_{(d)}(\mathbb{Z}[i])$ that belong to the ball of radius $H$ centered at the origin in $\mathcal{H}_{(d)}$ with respect to the unitarily invariant norm $\|\cdot\|_\Delta$. Namely,*

$$\overline{N}(1, H) := \sharp\left(B_\Delta(0, H) \bigcap \mathcal{H}_{(d)}(\mathbb{Z}[i])\right).$$

*Then, the following holds :*

$$|\overline{N}(1, H) - Vol(B_\Delta(0, H))| \leq \frac{\mathfrak{S}^{(2N+2)}}{\det(\Delta_{(d)})^2} H^{2N+1}.$$

*In particular, we conclude the following inequalities*

$$\frac{K_{2N+2}}{\det(\Delta_{(d)})} H^{2N+2} - \frac{\mathfrak{S}^{(2N+2)}}{\det(\Delta_{(d)})^2} H^{2N+1} \leq N(1, H).$$

## 5.2 On $\mathbb{C}$–visible ill–conditioned systems of polynomial equations

In this Subsection we use the notion of $\mathbb{C}$–visible points (Subsection 4.4 above) to count points in $\mathbb{P}(\mathcal{H}_{(d)}(\mathbb{Q}[i]))$.

**Proposition 37** *Let $H, \varepsilon \in \mathbb{R}$ be two positive real numbers, $H \geq 1$. Let $\overline{\mathcal{N}}(\varepsilon, H)$ be the number of points in $\mathbb{P}(\mathcal{H}_{(d)}(\mathbb{Q}[i]))$ of unitarily invariant height at most $H$ which belong to $\Sigma(\varepsilon)$. Then, the following inequality holds :*

$$|\overline{\mathcal{N}}(\varepsilon, H) - \frac{Vol(\overline{R}_\Delta(\varepsilon, H))}{4\zeta(2N+2)}| \leq \overline{L}(\varepsilon, N) H^{2N+1} + \frac{H}{4},$$

*where*

$$\overline{L}(\varepsilon, N) = \frac{\overline{T}_N \mathfrak{S}^{(2N+2)}}{4 \det(\Delta_{(d)})^2 \zeta(2N+1)} + \frac{Vol(\overline{R}_\Delta(\varepsilon, 1))}{2},$$

*and $\overline{T}_N$ is the constant defined in Proposition 35 above.*
*In particular, the following inequality is a consequence of Corollary 29:*

$$\overline{\mathcal{N}}(\varepsilon, H) \leq \frac{\varepsilon^4 \mathfrak{D}[(d)]}{4\zeta(2N+2)} H^{2N+2} + \overline{L}'(\varepsilon, N) H^{2N+1} + \frac{H}{4},$$

*where $\mathfrak{D}[(d)]$ is the constant introduced in Corollary 29 and*

$$\overline{L}'(\varepsilon, N) := \frac{\overline{T}_N \mathfrak{S}^{(2N+2)}}{4 \det(\Delta_{(d)})^2 \zeta(2N+1)} + \frac{\varepsilon^4 \mathfrak{D}[(d)]}{2}.$$

*Proof.–* Assume that $\varepsilon > 0$ is fixed throughout this proof. We also identify $\mathcal{H}_{(d)} \cong \mathbb{C}^{N+1}$ and $\mathcal{H}_{(d)}(\mathbb{Z}[i]) \cong \mathbb{Z}[i]^{N+1}$. Let $\overline{V}(\varepsilon) \subseteq \mathbb{C}^{N+1}$ be the compact subset given by the following identity :

$$\overline{V}(\varepsilon) := B_\Delta(0, 1) \bigcap \left(\pi_\mathbb{C}^{-1}(\Sigma(\varepsilon)) \cup \{0\}\right) = \overline{R}_\Delta(\varepsilon, 1).$$

For every positive real number $t \in \mathbb{R}_+$, the following identity holds :

$$t\overline{V}(\varepsilon) := \{tF : F \in \overline{V}(\varepsilon)\} \cong \overline{R}(\varepsilon, t).$$

In particular, for every positive real number $t \in \mathbb{R}_+$ the following identity also holds :

$$t^{2N+2} Vol(\overline{V}(\varepsilon)) = Vol(\overline{R}(\varepsilon, t)). \tag{27}$$

As in Subsection 4.4, for every positive real number $\rho \in \mathbb{R}_+$, let $\mathbb{Z}[i]^{N+1}(\rho)$ be the $\mathbb{Z}[i]$–lattice in $\mathbb{C}^{N+1}$ given by the following identity :

$$\mathbb{Z}[i]^{N+1}(\rho) := \{\rho X : X \in \mathbb{Z}[i]^{N+1}\}.$$

Let us define the finite set $A_\rho \subseteq B_\Delta(0, 1)$ given by the following equality

$$A_\rho := \overline{V}(\varepsilon) \bigcap \mathbb{Z}[i]^{N+1}(\rho) \setminus \{0\}.$$

Let $g(\rho)$ be the number of points in $A_\rho$. Observe that the following equality holds :

$$g(\rho) := \sharp(A_\rho) = \sharp\left(\overline{R}_\Delta(\varepsilon, \rho^{-1}) \bigcap \mathbb{Z}[i]^{N+1}\right) - 1 = \overline{N}(\varepsilon, \rho^{-1}) - 1. \tag{28}$$

For every $\lambda \in \mathbb{Z}[i]$, let us define the set $\Omega(\lambda, \rho)$ in the following terms :

$$\Omega(\lambda, \rho) := \{\lambda \rho X \in \overline{V}(\varepsilon) : gcd_{\mathbb{Z}[i]}(X) \in \mathbb{Z}[i]^*\},$$

where $gcd_{\mathbb{Z}[i]}(X)$ is the greatest common divisor of the coordinates of $X$ and $\mathbb{Z}[i]^* := \{1, -1, i, -i\}$ (cf. Subsection 4.4 above for more details).
For every $m \in \mathbb{N}$, $m \geq 1$, let us define :

$$\Omega(m, \rho) := \bigcup_{|\lambda|=\sqrt{m}} \Omega(\lambda, \rho).$$

The following holds :
Fact .– The following is a decomposition of $A_\rho$ as a disjoint union of subsets :

$$A_\rho = \bigcup_{m \geq 1} \Omega(m, \rho). \tag{29}$$

As $\Omega(m, \rho) \subseteq A_\rho$, this Fact is an immediate consequence of the following two claims :

i) Claim 1.– Given $m, m' \in \mathbb{N}$, $m \geq 1, m' \geq 1$, if $\Omega(m, \rho) \cap \Omega(m', \rho) \neq \emptyset$, then $m = m'$,

ii) Claim 2.– $A_\rho \subseteq \bigcup_{m \geq 1} \Omega(m, \rho)$.

In order to prove Claim 1, let $m, m' \in \mathbb{N}$ be such that $m \geq 1$, $m' \geq 1$ and $\Omega(m, \rho) \cap \Omega(m', \rho) \neq \emptyset$. Assume $Y \in \Omega(m, \rho) \cap \Omega(m', \rho)$. Then, there are $\lambda, \lambda' \in \mathbb{Z}[i]$ such that the following holds :

$$Y := \lambda \rho X_0 = \lambda' \rho X_0',$$

where $|\lambda| = \sqrt{m}$, $|\lambda'| = \sqrt{m'}$ and

$$gcd_{\mathbb{Z}[i]}(X_0) \in \mathbb{Z}[i]^*, \; gcd_{\mathbb{Z}[i]}(X_0') \in \mathbb{Z}[i]^*.$$

From Subsection 4.4, these last properties mean that $X_0$ and $X_0'$ are $\mathbb{C}$−visible from the origin. In particular,

$$\|X_0\|_\Delta = H_\Delta(X_0), \; \|X_0'\|_\Delta = H_\Delta(X_0').$$

Now, observe that $\pi_\mathbb{C}(Y) = \pi_\mathbb{C}(X_0) = \pi_\mathbb{C}(X_0')$. Thus, from Weil's product formula we conclude the following identity :

$$\|X_0\|_\Delta = H_\Delta(X_0) = H_\Delta(Y) = H_\Delta(X_0') = \|X_0'\|_\Delta.$$

Finally, the following holds :

$$\|Y\|_\Delta := |\lambda|\rho\|X_0\|_\Delta = |\lambda'|\rho\|X_0'\|_\Delta \Rightarrow \sqrt{m} = |\lambda| = |\lambda'| = \sqrt{m'},$$

and Claim 1 follows. As for Claim 2, let $\rho X \in A_\rho$ be a point in the lattice $\mathbb{Z}[i]^{N+1}(\rho)$. With the same notations as in Subsection 4.4, let $\lambda \in \mathbb{Z}[i]$ the Gauss integer given by the following identity :

$$\lambda := gcd_{\mathbb{Z}[i]}(X).$$

Then, there is $X_0 \in (\mathbb{Z}[i])^{N+1}$ such that $gcd_{\mathbb{Z}[i]}(X_0) \in \mathbb{Z}[i]^*$ and $X = \lambda X_0$. In particular,

$$\rho X = \lambda \rho X_0 \in \Omega(\lambda, \rho) \subseteq \Omega(|\lambda|, \rho).$$

Taking $m := |\lambda|^2 \geq 1$, Claim 2 follows.

For every integer number $m \in \mathbb{N}$, $m \geq 1$, let $f(m\rho)$ be the number of points in $\Omega(m, \rho)$. Observe that this function is well–defined for fixed $\rho$. Then, Equality (29) above, yields the following equality :

$$g(\rho) := \sum_{m=1}^{\infty} f(m\rho).$$

Applying Möbius inversion formula (see Theorem 270 in [14], for instance) it follows that

$$f(\rho) := \sum_{m=1}^{\infty} \mu(m) g(m\rho),$$

where $\mu$ is Möbius function. From Theorem 287 of [14] (cf. Equation (11) above) we have

$$\rho^{2N+2} f(\rho) - \frac{Vol(\overline{V}(\varepsilon))}{\zeta(2N+2)} = \sum_{m=1}^{\infty} \left( \rho^{2N+2} \mu(m) g(m\rho) - \frac{\mu(m)}{m^{2N+2}} \right).$$

Hence, we conclude

$$\rho^{2N+2} f(\rho) - \frac{Vol(\overline{V}(\varepsilon))}{\zeta(2N+2)} = \sum_{m=1}^{\infty} \frac{\mu(m)}{m^{2N+2}} \left( (m\rho)^{2N+2} g(m\rho) - Vol(\overline{V}(\varepsilon)) \right).$$
(30)

As in the Proof of Proposition 18, let $\mathbb{N}_\rho$ be the set of integer numbers given by the following identity

$$\mathbb{N}_\rho := \{m \in \mathbb{N} \;:\; 1 \leq m, m\rho \leq 1\}.$$

Let us denote by $S_\rho$ the sum of the terms of the left hand series of Equation (30) whose indices are in $\mathbb{N}_\rho$. Namely,

$$S_\rho := \sum_{m \in \mathbb{N}_\rho} \frac{\mu(m)}{m^{2N+2}} |(m\rho)^{2N+2} g(m\rho) - Vol(\overline{V}(\varepsilon))|$$

As in the Proof of Proposition 18 above, assume $m \in \mathbb{N}_\rho$. Then, replacing $(m\rho)^{-1}$ by $H$, we conclude from Identities (27) and (28) above the following identity :

$$|(m\rho)^{2N+2}g(m\rho) - Vol(\overline{V}(\varepsilon))| = \frac{1}{H^{2N+2}}|\overline{N}(\varepsilon, H) - (1 + Vol(\overline{R}_\Delta(\varepsilon, H)))| \tag{31}$$

From Proposition 35 above, we have :

$$\frac{1}{H^{2N+2}}|\overline{N}(\varepsilon, H) - Vol(\overline{R}_\Delta(\varepsilon, H))| \leq \frac{\overline{T}_N \mathfrak{S}^{(2N+2)}}{H \det(\Delta_{(d)})^2}.$$

Replacing back $H$ by $(m\rho)^{-1}$, Equation (31) becomes the following inequality :

$$|(m\rho)^{2N+2}g(m\rho) - Vol(\overline{V}(\varepsilon))| \leq (m\rho)\frac{\overline{T}_N \mathfrak{S}^{(2N+2)}}{\det(\Delta_{(d)})^2} + (m\rho)^{2N+2}. \tag{32}$$

We conclude :

$$S_\rho \leq \sum_{m \in \mathbb{N}_\rho} \left( \frac{\mu(m)}{m^{2N+1}} \rho \frac{\overline{T}_N \mathfrak{S}^{(2N+2)}}{\det(\Delta_{(d)})^2} + \mu(m)\rho^{2N+2} \right). \tag{33}$$

Then,

$$S_\rho \leq \frac{\rho \overline{T}_N \mathfrak{S}^{(2N+2)}}{\zeta(2N+1)\det(\Delta_{(d)})^2} + \sum_{m \in \mathbb{N}_\rho} \rho^{2N+2},$$

and we conclude

$$S_\rho \leq \frac{\rho \overline{T}_N \mathfrak{S}^{(2N+2)}}{\zeta(2N+1)\det(\Delta_{(d)})^2} + \rho^{2N+1}. \tag{34}$$

As in the Proof of Proposition 18, let $\mathbb{N}'_\rho$ be the class given by

$$\mathbb{N}'_\rho := \{m \in \mathbb{N} : 1 \leq m,\ m\rho > 1\}.$$

Observe that if $m \in \mathbb{N}'_\rho$, then $g(m\rho) = 0$.
Hence, let $S'_\rho$ be the sum of the terms in the left hand series of Equation (30) whose indices are in $\mathbb{N}'_\rho$. Namely,

$$S'_\rho := \sum_{m \in \mathbb{N}'_\rho} \frac{\mu(m)}{m^{2N+2}} |(m\rho)^{2N+2}g(m\rho) - Vol(\overline{V}(\varepsilon))|.$$

Then, we have

$$S'_\rho \leq \sum_{m \in \mathbb{N}'_\rho} \frac{\mu(m)}{m^{2N+2}} |Vol(\overline{V}(\varepsilon))| \leq \rho Vol(\overline{V}(\varepsilon)) \tag{35}$$

Now, we combine Equations (34) and (35) with Equation (30) to conclude

$$|\rho^{2N+2}f(\rho) - \frac{Vol(\overline{V}(\varepsilon))}{\zeta(2N+2)}| \leq \frac{\rho \overline{T}_N \mathfrak{S}^{(2N+2)}}{\zeta(2N+1)\det(\Delta_{(d)})^2} + \rho^{2N+1} + \rho Vol(\overline{V}(\varepsilon)).$$

Finally, replacing $\rho^{-1}$ by $H$ in this Equation and using the identity described in Equation 27 we conclude :

$$|f(H^{-1}) - \frac{Vol(\overline{R}_\Delta(\varepsilon, H))}{\zeta(2N+2)}| \leq \frac{\overline{T}_N \mathfrak{S}^{(2N+2)}}{\zeta(2N+1)\det(\Delta_{(d)})^2}H^{2N+1} + H + \frac{Vol(\overline{R}_\Delta(\varepsilon, H))}{H}.$$

Finally, from Lemma 33 above we have $\frac{f(H^{-1})}{4} = \overline{\mathcal{N}}(\varepsilon, H)$ and Proposition 37 is achieved. ■

The same arguments yield the following Proposition 19.

**Proposition 38** *Let $H, \varepsilon \in \mathbb{R}$ be two positive real numbers, $H \geq 1$. Let $\overline{\mathcal{N}}(1, H)$ be the number of points in $\mathbb{P}(\mathcal{H}_{(d)}(\mathbb{Q}[i]))$ of unitarily invariant height at most $H$. Then, the following inequality holds :*

$$|\overline{\mathcal{N}}(1, H) - \frac{Vol(B_\Delta(0, H))}{4\zeta(2N+2)}| \leq \overline{L}(1, N)H^{2N+1} + \frac{H}{4},$$

*where*

$$\overline{L}(1, N) = \frac{\mathfrak{S}^{(2N+2)}}{4\det(\Delta_{(d)})^2 \zeta(2N+1)} + \frac{K_{2N+2}}{2\det(\Delta_{(d)})},$$

*In particular, the following inequality also holds :*

$$\frac{K_{2N+2}}{4\det(\Delta_{(d)})\zeta(2N+2)}H^{2N+2} - \left(\overline{L}(1, n)H^{2N+1} + \frac{H}{4}\right) \leq \overline{\mathcal{N}}(1, H)$$

### 5.3 Proof of Theorem 4

We are now in conditions to prove Theorem 4 as stated at the Introduction.
*Proof of Theorem 4.–* First of all, observe that the probability that a random choice of a system of polynomial equations $F \in \mathbb{P}(\mathcal{H}_{(d)}\mathbb{Q}[i])$ of unitarily invariant height at most $H$ satisfies $\mu_{norm}(F) \geq 1/\varepsilon$ is given by the following quotient :

$$\frac{\overline{\mathcal{N}}(\varepsilon, H)}{\overline{\mathcal{N}}(1, H)}.$$

Let $\overline{A}$ be the quantity :

$$\overline{A} := \frac{K_{2N+2}}{4\det(\Delta_{(d)})\zeta(2N+2)}.$$

From Propositions 38 and 37 we conclude :

$$\overline{\mathcal{N}}(\varepsilon, H) \leq \overline{A}\mathfrak{C}[(d)]\varepsilon^4 H^{2N+2} + \overline{L}'(\varepsilon, N)H^{2N+1} + \frac{H}{4}.$$

$$\overline{\mathcal{N}}(1, H) \geq \overline{A}H^{2N+2} - \left(\overline{L}(1, N)H^{2N+1} + \frac{H}{4}\right),$$

where $\mathfrak{C}[(d)] := n^3(n+1)N(N-1)\mathcal{D}_{(d)}$ is the constant introduced in Theorem 3 of the Introduction. Then, we conclude

$$\frac{\overline{\mathcal{N}}(\varepsilon, H)}{\overline{\mathcal{N}}(1, H)} \leq \mathfrak{C}[(d)]\varepsilon^4 + \frac{\mathcal{F}(\varepsilon, N, H)}{H - \mathcal{G}(N, H)},$$

where

$$\mathcal{G}(N, H) := \frac{\overline{L}(1, N)}{\overline{A}} + \frac{1}{4\overline{A}H^{2N}},$$

and

$$\mathcal{F}(\varepsilon, N, H) := \varepsilon^4 \mathfrak{C}[(d)]\mathcal{G}(N, H) + \frac{\overline{L}'(\varepsilon, N)}{\overline{A}} + \frac{1}{4\overline{A}H^{2N}}.$$

Simplifying these expressions we conclude :

$$\mathcal{G}(N, H) = \zeta(2N+2)\left[\frac{\mathfrak{S}^{(2N+2)}}{\det(\Delta_{(d)})\zeta(2N+1)K_{2N+2}} + \frac{\det(\Delta_{(d)})}{K_{2N+2}H^{2N}} + 2\right].$$

On the other hand, the following identity holds :

$$\frac{\overline{L}'(\varepsilon, N)}{\overline{A}} = \zeta(2N+2)\left[\frac{\overline{T}_N\mathfrak{S}^{(2N+2)}}{\det(\Delta_{(d)})\zeta(2N+1)K_{2N+2}} + 2\varepsilon^4\mathfrak{C}[(d)]\right].$$

Hence, we conclude

$$\mathcal{F}(\varepsilon, N, H) = \zeta(2N+2)\left[\frac{(\overline{T}_N + \varepsilon^4\mathfrak{C}[(d)])\mathfrak{S}^{(2N+2)}}{\det(\Delta_{(d)}\zeta(2N+1)K_{2N+2}}\right] +$$

$$+ \zeta(2N+2)\left[\left(\varepsilon^4\mathfrak{C}[(d)] + 1\right)\frac{\det(\Delta_{(d)})}{K_{2N+2}H^{2N}} + 4\varepsilon^4\mathfrak{C}[(d)]\right].$$

As $\det(\Delta_{(d)}) \leq 1$ and $H \geq 1$, we may conclude :

$$\frac{\overline{\mathcal{N}}(\varepsilon, H)}{\overline{\mathcal{N}}(1, H)} \leq \mathfrak{C}[(d)]\varepsilon^4 + \frac{\mathcal{F}(\varepsilon, N)}{\zeta(2N+2)^{-1}H - \mathcal{G}(N)},$$

where

$$\mathcal{G}(N) := \frac{\mathfrak{S}^{(2N+2)}}{\det(\Delta_{(d)})\zeta(2N+1)K_{2N+2}} + \frac{1}{K_{2N+2}} + 2,$$

and

$$\mathcal{F}(\varepsilon, N) := \frac{(\overline{T}_N + \varepsilon^4\mathfrak{C}[(d)])\mathfrak{S}^{(2N+2)}}{\det(\Delta_{(d)}\zeta(2N+1)K_{2N+2}} + \frac{\left(\varepsilon^4\mathfrak{C}[(d)] + 1\right)}{K_{2N+2}} + 4\varepsilon^4\mathfrak{C}[(d)].$$

■

# 6  Applications : On the average precision required to write down an Approximate Zero

In this Section we discuss some applications of the estimates of Theorem 4 to compute the average precision required to write down approximate zeros for the affine Newton operator. We start by recalling some elementary facts about these ideas.

## 6.1  Some well–known facts about Affine Approximate Zero Theory.

Let $(d) := (d_1, \ldots, d_n)$ be a list of degrees, and let $\mathcal{P}_{(d)}$ be the set of all sequences $F := (f_1, \ldots, f_n)$ of multivariate polynomials that satisfy the following properties for every $i$, $1 \leq i \leq n$ :

i) $f_i \in \mathbb{C}[X_1, \ldots, X_n]$ and

ii) $\deg(f_i) \leq d_i$,

We may easily identify $\mathcal{P}_{(d)}$ and $\mathcal{H}_{(d)}$ by the obvious homogenization operator. For every $F \in \mathcal{P}_{(d)}$, we also denote by $F$ its homogenization $F \in \mathcal{H}_{(d)}$. In this Subsection we study the affine counterpart of multivariate Newton's operator. Let $F := (f_1, \ldots, f_n) \in \mathcal{P}_{(d)}$ be a sequence of multivariate polynomials. Let $V(F) \subseteq \mathbb{C}^n$ be the set of common zeros into the complex affine space $\mathbb{C}^n$, namely

$$V(F) := \{\zeta \in \mathbb{C}^n \ : \ f_1(\zeta) = \cdots = f_n(\zeta) = 0\}.$$

Let $\zeta \in V(F)$ be a smooth zero of system $F$, namely the jacobian matrix of $F$ at $\zeta$ is non–singular, i.e. $DF(\zeta) \in GL(n, \mathbb{C})$, where

$$DF(\zeta) := \left( \frac{\partial f_i}{\partial X_j}(\zeta) \right)_{1 \leq i,j \leq n}.$$

Let $N_F$ be the multivariate Newton operator defined by $F$. Namely,

$$N_F(X_1, \ldots, X_n) := \begin{pmatrix} X_1 \\ \vdots \\ X_n \end{pmatrix} - DF(X_1, \ldots, X_n)^{-1} \begin{pmatrix} f_1(X_1, \ldots, X_n) \\ \vdots \\ f_n(X_1, \ldots, X_n) \end{pmatrix},$$

Next, we recall Smale's Approximate Zero Theory (cf. [33], [28], [27], [29], [30], [32], [31] and the compiled version in [3]) as foundation of Numerical Analysis.

With the same notations and assumptions as above, an *approximate zero* of the system $F$ with associated zero $\zeta$ is a point $z \in \mathbb{C}^n$ such that the sequence of iterates of the Newton operator is well–defined and converges quadratically to $\zeta$. Namely, an approximate zero of the system $F$ with associated zero $\zeta \in V(F)$ is an affine point $z \in \mathbb{C}^n$ such that the following properties hold :

- The jacobian matrix of $F$ at $z$ is a regular matrix. Namely, $DF(z) \in GL(n, \mathbb{C})$.

- The following sequence is well–defined :

$$z_1 := N_F(z) \in \mathbb{C}^n, \text{ and } z_k := N_F(z_{k-1}) \text{ for } k \geq 2.$$

- For every $k \in \mathbb{N}$, $k \geq 1$, the following inequality holds :

$$\|z_k - \zeta\| \leq \frac{1}{2^{2^{k-1}}} \|z - \zeta\|$$

With the same notations as above, we define the *quantity* $\gamma$ in the following terms :

$$\gamma(F, \zeta) := \sup_{k \geq 2} \left\| \frac{(DF(\zeta))^{-1}(D^{(k)}F(\zeta))}{k!} \right\|^{\frac{1}{k-1}},$$

where the norm on the right hand side of this identity is the norm of the multilinear operator

$$DF(\zeta)^{-1} D^{(k)} F(\zeta) \ : (\mathbb{C}^n)^k \longrightarrow \mathbb{C}^n.$$

This quantity yields a locally sufficient condition for having an approximate zero.

**Theorem 39 ($\gamma$–Theorem, [33])** *With the same notations as above, Let $F \in \mathcal{P}_{(d)}$ be a regular sequence of polynomials that defines a smooth, zero–dimensional affine algebraic variety $V(F)$, and let $\zeta \in \mathbb{C}^n$ a smooth zero of system $F$. Let $z \in \mathbb{C}^n$ be an affine point satisfying the following inequality :*

$$\|\zeta - z\| \gamma(f, \zeta) \leq \frac{3 - \sqrt{7}}{2}.$$

*Then, $z$ is an approximate zero of the system $F$ with associated zero $\zeta$.*

For every affine point $\zeta := (\zeta_1, \ldots, \zeta_n) \in \mathbb{C}^n$, let us denote by $\widetilde{\zeta} \in \mathbb{P}_n(\mathbb{C})$ the corresponding projective point given by the following identity :

$$\widetilde{\zeta} := (1 : \zeta_1 : \cdots : \zeta_n) \in \mathbb{P}_n(\mathbb{C}).$$

Combining Lemma 7 and Theorem 2 of Chapter 14 in [3], for every $\zeta \in \mathbb{C}^n$, such that $F(\zeta) = 0$ and $DF(\zeta) \in GL(n, \mathbb{C})$, the following inequality holds :

$$\gamma(F, \zeta) \leq \frac{(D_{(d)})^{3/2} \mu_{norm}(F, \widetilde{\zeta})}{2}. \tag{36}$$

## 6.2 On the number of approximate zeros of given precision

Let us assume now the *computational hypothesis* (cf. also [4]) : Assume your are given a system of polynomial equations $F \in \mathcal{P}_{(d)}$ and you want to compute an approximate zero $z \in \mathbb{C}^n$ of system $F$ with associated zero $\zeta \in \mathbb{C}^n$. From a computational point of view, you look for approximate zeros $z$ that could be written by a computer (or equivalently that can be written as a list of symbols in a finite alphabet). In standard numerical analysis procedures, the computational hypothesis becomes now the assumption that $z \in \mathbb{Q}[i]^n$.

Given a point $z \in \mathbb{Q}[i]^n$, we define *the precision of $z$* as the bit length of its denominator. Namely, let $z \in \mathbb{Q}[i]^n$ be a point whose coordinates are Gauss rationals. Then, there are $q \in \mathbb{Z} \setminus \{0\}$ and $(z_1, \ldots, z_n) \in \mathbb{Z}[i]^n$ such that

$$z := \left(\frac{z_1}{q}, \ldots, \frac{z_n}{q}\right),$$

and $q$ is of minimal absolute value satisfying this property. The precision $Pr(z)$ of $z \in \mathbb{Q}[i]^n$ is defined to be $Pr(z) := \max\{0, \log_2 q\}$.

The following statement gives upper and lower bounds for the number of approximate zeros of given precision that satisfy the $\gamma-$Theorem above. In other words, we show upper and lower estimates for the behaviour of the function :

$$N_m(\zeta) := \sharp\{z \in \mathbb{Q}[i]^n : \|z - \zeta\| \leq \frac{3 - \sqrt{7}}{2\gamma(F, \zeta)}, \wedge \, Pr(z) = \log \, m\}.$$

Observe that $N_m(\zeta)$ equals the cardinal of the set :

$$B(\zeta, \frac{(3 - \sqrt{7})m}{\gamma(F, \zeta)}) \bigcap \mathbb{Z}[i]^n.$$

**Proposition 40** *With the same notations as above, the following properties hold :*

*i) For every integer number $H$ such that*

$$H < \left(\frac{\gamma(F, \zeta)}{3 - \sqrt{7}}\right)^{\frac{1}{2}},$$

*the following holds*

$$0 \leq \sum_{m=1}^{H} N_m(\zeta) \leq 1.$$

*ii)* For every $m \in \mathbb{N}$ such that
$$H_1 := \left(\frac{\gamma(F,\zeta)}{3-\sqrt{7}}\right)^{\frac{1}{2}} \leq m \leq H_2 := \frac{\gamma(F,\zeta)}{3-\sqrt{7}},$$
*the following holds*
$$0 \leq N_m(\zeta) \leq 1.$$
*In particular,*
$$0 \leq \sum_{m=H_1}^{H_2} N_m(\zeta) \leq H_2 - H_1 = \frac{\gamma(F,\zeta)}{3-\sqrt{7}} - \sqrt{\frac{\gamma(F,\zeta)}{3-\sqrt{7}}}.$$

*iii)* *For every $m \in \mathbb{N}$ the following inequalities hold :*

$$K_n \left(\frac{m(3-\sqrt{7})}{2\gamma(F,\zeta)} - \sqrt{2n}\right)^{2n} \leq N_m(\zeta) \leq K_n \left(\frac{m(3-\sqrt{7})}{2\gamma(F,\zeta)} + \frac{\sqrt{2n}}{2}\right)^{2n}$$

*Proof.–* Item *iii)* follows from an elementary argument as those of [20] or the more classical Blichfeldt–Minkowski estimates. Similar estimates could also be obtained by means of Davenport's Theorem (Theorem 12 above). We left the proof of item *iii)* for the reader. As for items *i)* and *ii)* we make use of the following estimate, which is a translation of the "Gap Principle" of [21] to our context (cf. also [26, 25]).
Given $z, z' \in \mathbb{Q}[i]^n$ two different approximate zeros of system $F$ with associate zero $\zeta$ that satisfy the $\gamma$–Theorem, and such that $mz, mz' \in \mathbb{Z}[i]^n$, the following chain of inequalities hold

$$\frac{1}{m} \leq \|\frac{mz}{m} - \frac{mz'}{m}\| \leq \|z - \zeta\| + \|z' - \zeta\| \leq \frac{3-\sqrt{7}}{\gamma(F,\zeta)}.$$

Hence, we conclude $m \geq \frac{\gamma(F,\zeta)}{3-\sqrt{7}}$ and the following holds for every $m$ such that $1 \leq m < \frac{\gamma(F,\zeta)}{3-\sqrt{7}}$ :
$$N_m(\zeta) \leq 1. \tag{37}$$
Item *i)* follows from the following Remark. Assume that $H$ satisfies the following inequality
$$H < \left(\frac{\gamma(F,\zeta)}{3-\sqrt{7}}\right)^{\frac{1}{2}}.$$
Let $z, z' \in \mathbb{Q}[i]^n$ be two points such that there are $d, d' \in \mathbb{N} \setminus \{0\}$ satisfying the following properties :

- The coordinates of the points $dz, d'z' \in \mathbb{Z}[i]^n$ are Gauss integers,

- $z, z'$ are approximate zeros of system $F$ with associate zero $\zeta$ that satisfy the $\gamma$−Theorem above. Namely, assume the following holds

$$\|z - \zeta\| \leq \frac{3 - \sqrt{7}}{2\gamma(F,\zeta)}, \quad \|z' - \zeta\| \leq \frac{3 - \sqrt{7}}{2\gamma(F,\zeta)},$$

- $d, d' < H$.

Then, the following chain of inequalities hold :

$$\frac{1}{dd'} \leq \|z - z'\| \leq \|z - \zeta\| + \|z' - \zeta\| \leq \frac{3 - \sqrt{7}}{\gamma(F,\zeta)}.$$

In particular, we would conclude $dd' \geq \frac{\gamma(F,\zeta)}{3-\sqrt{7}}$. In other words, if $d \leq H$, then $d' \geq H$ and item $i$) follows. Item $ii$) follows from the previous Remark and Inequality 37 above. ∎

Item $iii$) of Proposition 40 immediately yields the following

**Corollary 41** *With the same notations and assumptions as above, for every $p \in \mathbb{N}$ such that the following inequality holds :*

$$p \geq \log \gamma(F,\zeta) + \log\left[K_n^{-1/2n} + \sqrt{2n}\right] + 1,$$

*there are approximate zeros $z \in \mathbb{Q}[i]^n$ of system $F$ of precision $p$ with associated zero $\zeta$.*

The following statement follows from the previous Corollary, Identity 36 and our estimates on the probability distribution of $\mu_{norm}$.

**Corollary 42** *With the same notations as above, there is a universal constant $c_3 > 0$ ( $c_3 \leq 20$) such that the following holds :*
*Let $(d) := (d_1, \ldots, d_n)$ be a list of degrees and let $h, w \in \mathbb{R}$ be two positive real numbers. Assume that $h$ is big enough to satisfy the following inequality :*

$$h \geq c_3 N^2 (\log N + \log D_{(d)}) + \log w,$$

*Then, for every system of polynomial equations $F \in \mathbb{P}(\mathcal{P}_{(d)}(\mathbb{Q}[i]))$ of bit length at most $h$, the required precision $Pr(F)$ to write down an approximate zero of $F$ satisfies*

$$Pr(F) \leq O(n \log_2 D_{(d)} + \log_2 w), \tag{38}$$

*with probability at least*

$$1 - \frac{2}{w}.$$